\documentclass[a4paper,12pt]{amsart} 

\usepackage{soul}
\usepackage{amsmath, amssymb, amsfonts, amsthm, amssymb}
\usepackage{euscript, mathrsfs, stackrel}
\usepackage{enumerate}
\usepackage{stmaryrd}
\usepackage[utf8]{inputenc}   
\usepackage[T1]{fontenc}      
\usepackage{ulem} 

\usepackage{color}
\usepackage{graphicx}
\usepackage{hyperref}
\usepackage{wrapfig}
\usepackage{cleveref}

\theoremstyle{plain}

\numberwithin{equation}{section}

\newtheorem{theorem}{Theorem}[section]

\newtheorem{lemma}[theorem]{Lemma}
\newtheorem{fact}[theorem]{Fact}
\newtheorem{corollary}[theorem]{Corollary}
\newtheorem{remark}[theorem]{Remark}

\usepackage[left=2.5cm,right=2.5cm,vmargin=2.5cm]{geometry}

\usepackage{mathrsfs} 

\def\z{\mathbb Z}
\def\r{\mathbb R}
\def\e{\mathbb E} 
\def\p{\mathbb P} 

\def\L{{\mathscr L}} 

\def\T{{\mathcal T}} 

\def\dd{{\rm d}}

\def\cS{{\mathcal S}}

 \title{Boundedness of discounted tree sums}

 \date{This version:  \today}
 
\author{Elie A\"id\'ekon}
\address{Elie A\"id\'ekon, SMS, Fudan University, China}
\email{aidekon@fudan.edu.cn}

\author{Yueyun Hu}
\address{Yueyun Hu, 
LAGA, Universit\'e Paris XIII,  93430 Villetaneuse,
France}
\email{yueyun@math.univ-paris13.fr}

\author{Zhan Shi}
\address{Zhan Shi,  AMSS,
Chinese Academy of Sciences,  China}
\email{shizhan@amss.ac.cn}

\begin{document}

\begin{abstract}  Let $(V(u),\, u\in \T)$ be a (supercritical) branching random walk and $(\eta_u,\,u\in \T)$ be marks on the vertices of the tree,  distributed in an i.i.d.\ fashion. Following Aldous and Bandyopadhyay \cite{AB05}, for each infinite ray $\xi$ of the tree, we associate  the {\it discounted tree sum}  $D(\xi)$ which is the sum of the $e^{-V(u)}\eta_u$ taken along the ray. The paper deals with the finiteness of $\sup_\xi D(\xi)$. To this end, we study the extreme behaviour of the local time processes  of the paths $(V(u),\,u\in \xi)$. It answers a question of Nicolas Curien, and partially solves Open Problem 31 of Aldous and Bandyopadhyay \cite{AB05}. We also present several open questions.
\end{abstract}



\maketitle

\section{Introduction}

Let  $V:=\{ V(u), u \in \T\}$ be   a discrete-time  branching random walk on the real line $\r$, where $\T$ is a Ulam-Harris tree which describes the genealogy of the particles and $V(u)\in \r$ is the position of the particle $u$. When  a particle $u$ is at $n$-th generation, we write $\vert u \vert=n$ for   $n \ge0$. The branching random walk $V$ is constructed as follows:  initially, there is a single particle $\varnothing$ located at $0$. The particle $\varnothing$ is  considered as  the root of $\T$.  At the first generation, the root dies and  gives birth to a certain number of children, these children along with their spatial positions form a point process $\aleph$ on $\r$ and constitute the first generation of the branching random walk $\{ V(u), \vert u \vert=1\}$. We may identify $\aleph=\sum_{\vert u \vert=1} \delta_{\{V(u)\}}$. For the next generations, the process is constructed recursively:  for each $\vert u \vert=n$ (if such $u$ exists) with $n \ge 1$,   the  particle $u$  dies  in the  $(n+1)$-th generation   and  gives birth to    an independent copy of $\aleph$ shifted by $V(u)$.   The collection of all children of all $u$ together with their positions  gives the $(n+1)$-th generation.  The whole system may survive forever or die out after some generations.

\bigskip

Define the measure $\mu(\dd x):=\e[\aleph(\dd x)]$. We assume that $\mu(\r)>1$, so that $\T$ is a supercritical Galton--Watson tree. We denote by $\cS$ the event of non-extinction of $\T$, and $\p^*:=\p(\cdot \mid \cS)$. By a ray $\xi$ of $\T$, we mean a path $\xi$ that is either finite—namely, the shortest path from the root $\varnothing$ to a leaf of $\T$—or infinite, in which case $\xi=\{\xi_0=\varnothing, \xi_1, ..., \xi_n, ...\}$, with $\xi_n$ being the parent of $\xi_{n+1}$ for all $n \ge 0$.
Denote by $\partial\T$  the set of all rays (finite or infinite) of $\T$. 
Let 
\begin{equation}\label{def:M}
M:= \sup_{\xi\in \partial \T}  \sum_{n=0}^\infty e^{-V(\xi_n)}.
\end{equation}

\noindent The starting point of this work is the following question:   when do we have  $M<\infty,$ $\p^*$-a.s? This question was raised by Nicolas Curien (personal communication) in relation to the construction of self-similar Markov trees \cite{BCR}, see in particular the discussion in Section 3.4 there. 

\bigskip

Let $m_n:= \inf_{\vert u \vert=n} V(u)$ and $\varphi(t):=\log \int_\r e^{-t x} \mu(\dd x)$ for $t\in \r$. We have $M\ge \sup_{n\ge 0} e^{-m_n}$. Hence $M=\infty$ if $\liminf_{n\to \infty} m_n=-\infty$. This happens $\p^*$-a.s.\ when $\inf_{t>0} {\varphi(t)\over t} \in (0,\infty)$. Actually, as soon as $\inf_{t>0} {\varphi(t)\over t}<\infty$, we have the following law of large numbers (\cite{H74, K75, biggins-mart-cvg}): $\p^*$-a.s.,
\begin{equation}\label{LLN}
\lim_{n\to\infty} \frac{m_n}{n}= - \inf_{t>0}\frac{\varphi(t)}{t}=:\gamma.
\end{equation}

\noindent  Notice that the upper bound
\begin{equation}\label{eq:upper_bound}
 M \le \sum_{n=0}^\infty e^{-m_n}
\end{equation}

\noindent implies that $M <\infty$ a.s. when $\gamma>0$. The critical case is therefore $\gamma=0$, where the upper bound \eqref{eq:upper_bound}  is not good enough in general. From now on, we suppose that one of the following assumptions holds: \\

{\bf (H1)} There exists $t>0$ such that $\varphi(t)=0$ and we let $t^*$ be the minimal such $t$. The integral $\int_\r xe^{-t^*x}\mu(\dd x)$ is well-defined and nonzero.  \\

{\bf (H2)} There exists $t^*>0$ such that $\varphi(t^*)=0$  and $\varphi(t)\in [0,\infty]$ for all $t>0$. The integral $\int_\r xe^{-t^*x}\mu(\dd x)$ is well-defined, is equal to $0$ and  $\int_\r x^2 e^{-t^*x}\mu(\dd x)\in (0,\infty)$. \\

In  (H1)  and  (H2), we say that $\int_\r xe^{-t^*x}\mu(\dd x)$ is well-defined if at least one of the quantities $\int_{(-\infty,0)} |x|e^{-t^*x}\mu(\dd x)$ or $\int_{(0,\infty)} xe^{-t^*x}\mu(\dd x)$ is finite. Note that $\varphi(t^*)=0$ yields that $m_n \to \infty$, $\p^*$-a.s. (see \cite[Lemma 3.1]{Shi12})  and $\gamma\ge 0$.  

Let us discuss the case  (H1): if $\gamma>0$, then $\int_\r xe^{-t^*x}\mu(\dd x)>0$. Otherwise, one necessarily has $\gamma=0$.  If furthermore $\int_\r xe^{-t^*x}\mu(\dd x)$ is negative, resp. positive, then one necessarily has $\varphi(t)=\infty$ for all $t< t^*$, resp.\ for all $t>t^*$.   The case (H1)  with $\int_\r xe^{-t^*x}\mu(\dd x)>0$ and $\mu((-\infty,0))>0$ will be denoted by {\bf (H1')}. 
 
 In case  (H2), one necessarily has $\gamma=0$.  The case (H2) is called ``boundary case'' in the literature, see \cite{A13} for the precise asymptotic of $m_n$ and \cite{Shi12} for further references. There exist cases where $\gamma=0$ but neither   (H1) nor  (H2) holds, for example when $\mu((-\infty,0))=0$ and $\mu(\{0\})  >   1$. In that case $\gamma=0$ while there is no $t>0$ such that $\varphi(t)=0$.
 We refer to the appendix of the arXiv version of \cite{J12} for such discussions.

\begin{theorem}\label{t:main}
Under {\bf (H1)} or {\bf (H2)}, $M<\infty$ \mbox{\rm $\p^*$-a.s.} Actually, 
$$
\Xi_b:=\sup_{\xi\in\partial\T} \sum_{n=0}^\infty (1+|V(\xi_n)|)^{-b}<\infty
$$ 

\noindent for all $b>2$ under {\bf (H1)}, and for all $b>3$ under {\bf (H2)}.
\end{theorem}

\bigskip

In the case $\gamma>0$, $\Xi_b<\infty$ for $b>1$ by \eqref{LLN} and \eqref{eq:upper_bound}. The fact that there exist rays such that $\limsup_{n\to\infty} \frac{V(\xi_n)}{n}<\infty$ implies that $\Xi_b=\infty$ for $b=1$.  On the other hand, according to Jaffuel \cite{J12}, under (H2) with some further integrability conditions on $\mu$, there exists an explicit constant $a_0>0$ such that  $\p^*$-a.s.
$$
\inf_{\xi\in \partial \T} \limsup_{n\to\infty} \frac{V(\xi_n)}{n^{1/3}}=a_0.
$$

\noindent  It implies that $\Xi_b = \infty$ for $b=3$ in this case. Finally, the following theorem shows there are cases under (H1') where $\Xi_b =\infty$ for $b=2$.

\begin{theorem}\label{t:consistent}
We suppose that {\bf (H1')} holds and there exist $1<\alpha_1 \le \alpha_2$  such that $$\e\Big[\sum_{|u|=1} e^{-t^*V(u)} |V(u)|^{\alpha_1} \Big]< \infty$$ and for all $y>0$ large enough,
$$
 \p\left(\exists \, |u|=1\,:\, V(u)\le -y\right)\ge e^{-t^* y} y^{-\alpha_2} .
$$
We further assume that $$\e \left[\Sigma_V \log_+ \Sigma_V\right]< \infty,$$  where $\Sigma_V:= \sum_{|u|=1} e^{-t^*V(u)}$ and $\log_+ x:= \max (\log x, 0)$ for $x>0$. Then there exists a constant $a_1\in (0,\infty)$ such that $\p^*$-a.s. 
$$
 \inf_{\xi\in \partial \T} \limsup_{n\to\infty} \frac{V(\xi_n)}{\sqrt{n\log n}} =a_1.
$$
\end{theorem}

\bigskip

The issue of the finiteness of $M$ is a particular case of the following problem stated in Aldous and Bandyopadhyay  \cite{AB05} as Open Problem 31, Section 4.4. It is related to the recursive distributional equation 
\begin{equation}\label{RDE}
X \overset{d}{=} \eta + \max_{\vert u \vert=1} e^{-V(u)}X_u 
\end{equation}

\noindent where on the right-hand side, conditionally on $(\eta, V(u),\, \vert u \vert=1)$, the r.v. $(X_u,\, \vert u \vert=1)$ are i.i.d.\ with the law of $X$.  We refer to \cite[Section 4.5]{AB05} for various examples where the equation \eqref{RDE} appears.  

 A possible endogenous solution of \eqref{RDE} can be constructed as follows. Consider a pair $(\eta,\aleph)$ where $\eta>0$ is a positive r.v. and $\aleph$ is a point process as above.  We do not suppose that $\eta$ and $\aleph$ are independent. We then construct recursively the branching random walk as before: for each vertex $u$ at position $V(u)$, we take an independent copy of $(\eta,\aleph)$, call it $(\eta_u,\aleph_u)$, so that the children of $u$ are at positions $V(u)+\aleph_u$.  Following \cite{AB05},   we introduce  the {\it discounted tree sum} \begin{equation} \label{def-D}    D(\xi):= \sum_{n=0}^\infty e^{-V(\xi_n)}\eta_{\xi_n}, \qquad \forall\, \xi\in \partial\T. \end{equation}

\noindent Define 
$$
X:= \sup_{\xi\in \partial \T}  D(\xi).
$$

\noindent Then $X$ is a solution of \eqref{RDE}.  Indeed, $\sup_{\xi\in \partial \T}  D(\xi)$ is the minimal positive solution of \eqref{RDE}. Therefore the existence of a positive solution of \eqref{RDE} is equivalent to the boundedness of 
the discounted tree sums.

 Observe that $X=M$ if $\eta\equiv 1$. Theorem 32 of \cite{AB05} shows (among other results) that $X<\infty$ a.s.  when $\mu((-\infty,0])=0$,  $E[\eta^p]<\infty$ for all $p\ge 1$ and $\varphi(t)<\infty$ for some $t\ge 1$. It was left open to study general conditions 
under which $X<\infty$. Again, it is clear from \eqref{LLN} that $X=\infty$ $\p^*$-a.s. if $\gamma<0$. Notice that $\gamma>0$ no longer ensures that $X<\infty$ a.s. because of the influence of the variable $\eta$. The property $\{X<\infty\}$  is inherited, namely, $\{X< \infty\}= \cap_{\vert u \vert=1} \{X^{(u)} < \infty\}$, where  $X^{(u)}$ is defined as $X$ but for the branching random walk indexed by the subtree of $\T$ rooted at $u$.  Therefore (see \cite{LPbook}) $\p^*(X=\infty) \in \{0, 1\}$.

\begin{theorem}\label{t:discounted}
Suppose that 
\begin{equation}\label{taileta}
\zeta:=\lim_{x\to\infty} \frac{-\log \p(\eta>x)}{\log x}
\end{equation}

\noindent  exists and  $\zeta\in (0, \infty]$. Under {\bf (H1)} or {\bf (H2)}, 
\begin{enumerate}[(i)]
\item if $\zeta<t^*$, then $\p^*(X=\infty)=1$.
\item if $\zeta>t^*$, then $\p^*(X=\infty)=0$.
\end{enumerate}
\end{theorem}

\bigskip
  Note that we allow $\zeta=\infty$ in the assumption \eqref{taileta}. This occurs when the tail of $\eta$ decays faster than any polynomial rate, such as  when $\eta$ possesses finite exponential moments.  

\begin{remark}   \begin{enumerate}
\item Theorem \ref{t:discounted}~(i) comes from the fact that $\sup_{u\in \T} e^{-V(u)}\eta_u=\infty,$ $\p^*$-a.s., see Lemma \ref{l:supdiscounted}. We refer to \cite{BG23}, \cite{GM24} for much more precise results on this model, which is called last progeny modified branching random walk there. 
\item Theorem \ref{t:discounted}~(ii) was proved in \cite{athreya} when $V(u)= c |u|$ for some constant $c>0$.
\item If $\partial \T$ is equipped with the usual ultra metric, then in the case (ii), $\xi \mapsto D(\xi)$ is continuous on $\partial\T$. For more details, see Remark \ref{r:continuity}. 
\end{enumerate}
\end{remark}

\bigskip
The way we address these problems is by looking at the local times of the random walk along a ray defined as  for any $\xi\in \partial\T$ and $k\in \z$, 
\begin{equation}\label{def:localtime}
N_\xi^k:= \sum_{n=0}^\infty 1_{\{k\le V(\xi_n)<k+1 \}}.
\end{equation}

\noindent Theorem \ref{t:main} is then a consequence of the following result.
\begin{theorem}\label{t:localtime}
We have \mbox{\rm $\p^*$-a.s.} 
$$
\lim_{n\to \infty}\frac1{n^{\kappa}}\sup_{\xi\in \partial T} N_\xi^n=0
$$ 

\noindent for all $\kappa>1$ under {\bf (H1)} and  all $\kappa>2$ under {\bf (H2)}.
\end{theorem}

We can readily prove Theorem \ref{t:main}. \\

{\noindent\it Proof of Theorem \ref{t:main}}. Let $b>0$. As mentioned before, $\lim_{n\to\infty} m_n\to \infty$ $\p^*$-a.s under our assumptions. 
Since $\sum_{n=0}^\infty {\bf 1}_{\{V(\xi_n)\ge 0\}}(1+V(\xi_n))^{-b} \le \sum_{n=0}^\infty (1+n)^{-b} N_{\xi}^n$,   we may conclude by Theorem \ref{t:localtime}. $\Box$

\bigskip

\bigskip

To study the critical cases of Theorem \ref{t:localtime}, we will restrict our attention  to some integer-valued right-continuous branching random walks. Recall that (H1') is the assumption (H1) together with the positivity of $\int_{\r} xe^{-t^*x}\mu(\dd x)$ and $\mu((-\infty,0))>0$.  
\begin{theorem}\label{t:limsup}
We suppose that $\mu$ is supported on $\mathbb{Z}_-\cup\{1\}$ and
$$
\e\Big[ \sum_{|u|,|v|=1,u\neq v} e^{- t^*V(u)}e^{- t^* V(v)} \Big]<\infty.
$$

\noindent  Notice that $N_{\xi}^k$ defined in \eqref{def:localtime} is now simply equal to $ \sum_{n=0}^\infty 1_{\{V(\xi_n) =k\}}$ when $k\ge 0$. \\
(i) Under {\bf (H1')}, $$\lim_{n\to \infty}\frac1{n}\sup_{\xi\in \partial T} N_{\xi}^n=-\frac{t^*}{\log q}, \qquad \mbox{\rm $\p^*$-a.s.}, $$
\noindent where $q\in (0,1)$ is given by \eqref{def:q}.\\
(ii) Under {\bf (H2)},  then $$\lim_{n\to \infty}\frac1{n^2}\sup_{\xi\in \partial T} N_{\xi}^n= \frac{t^*}{2\theta} , \qquad \mbox{\rm $\p^*$-a.s.}, $$ 
\noindent where $\theta$ is given by \eqref{theta}.
\end{theorem}

  It implies that, for such branching random walks, $$\sup_{\xi\in\partial \T} \limsup_{n\to\infty} \frac1n N_{\xi}^n = -\frac{t^*}{\log q}, \qquad \mbox{ under (H1'),}$$  and $$\sup_{\xi\in\partial \T} \limsup_{n\to\infty} \frac1{n^2}N_{\xi}^n=\frac{t^*}{2\theta}, \qquad \mbox{ under (H2)}. $$ 
 
 \noindent One can ask the same question for $\liminf$ instead of $\limsup$.  Theorem \ref{t:consistent} suggests that the renormalisation should be different under its conditions.

\bigskip

{\bf Open question 1}.

(i)  Under the assumptions of Theorem \ref{t:consistent}, is $\sup_{\xi\in\partial \T} \liminf_{n\to\infty} \frac{\log n}{n} N_{\xi}^n \in (0,\infty)?$ What is its value?

 (ii) Under the assumptions of Theorem \ref{t:limsup} (ii), is $\sup_{\xi\in\partial \T} \liminf_{n\to\infty} \frac1{n^2} N_{\xi}^n \in (0,\infty)?$ What is its value?

\bigskip 
{\bf Open question 2}.  In Theorem \ref{t:consistent}, what is the value of  $a_1$?

\bigskip

{\bf Open question 3}. Let $a>0$. Under the conditions of Theorem \ref{t:limsup} (ii), what is the Hausdorff dimension of the rays $\xi \in \partial \T$ such that $\limsup_{n\to\infty} \frac1{n^2} N_{\xi}^n=a$? Same question with $\liminf$ instead of $\limsup$.

\bigskip

The rest of the paper is organized as follows: 
Section \ref{s:manytoone} introduces the usual many-to-one formula. Theorem \ref{t:localtime}, Theorem \ref{t:discounted},  and  Theorem \ref{t:limsup} are proved  in Section \ref{s:localtime}, Section \ref{s:discounted} and Section \ref{s:limsup},  respectively. A common tool in these proofs is the construction of some suitable optional lines in the sense of  \cite{BK04} by using the local times of the branching random walk.  The proof of Theorem \ref{t:consistent} is given in Section \ref{s:consistent}  and relies on a certain inhomogeneous Galton--Watson process, building on the fact that under (H1') and some additional assumptions, the minimum of the branching random walk is achieved by a large drop (see \cite{BHM}).  The Appendix \ref{s:est} contains some estimates on random walks which are used in the course of the proof of Theorem \ref{t:limsup}.

  Throughout this paper, $c, c', c''$, eventually with some subscripts, denote positive constants whose values may vary from one line to another.

\bigskip

{\bf Acknowledgments}. We thank Nicolas Curien for introducing us to this question and explaining the link with self-similar Markov trees. 

\bigskip

{\it Shortly before completion of this work, we learnt that Bastien Mallein (personal communication) proved that $M<\infty$ in the case of the branching Brownian motion by a similar argument, hence solving the question of Nicolas Curien in that case.}

\section{Many-to-one formula}
\label{s:manytoone}

Since $\varphi(t^*)=0$, $M_n:=\sum_{|u|=n} e^{-t^* V(u)}$ defines a martingale. For $u\in \T$ and $0\le k\le |u|$, we let $u_k$ be the ancestor of $u$ at generation $k$. In particular $u_0=\varnothing$. We recall a very convenient tool in the study of branching random walks, namely  the many-to-one formula, see \cite[Theorem 1.1]{Shi12}:

\begin{fact}[The many-to-one formula] Under (H1) or (H2), there exists a one-dimensional random walk $(S_k)_{k\ge 0}$ such that $S_0=0$ and for any $k\ge 1$ and any measurable function $f: \r^k \to\r_+$, we have \begin{equation}     \e \Big[\sum_{|u|=k} f(V(u_1), ..., V(u_k))\Big]= \e \left[ e^{t^* S_k} f(S_1, ..., S_k)\right], \label{many-to-one} \end{equation}
where the step distribution of $S$ is given by $$\e [h(S_1)]=  \e\Big[\sum_{|u|=1} e^{-t^* V(u)} h(V(u))\Big],$$ for any Borel bounded function $h$.  In particular, $$\e[S_1]=  \e\Big[\sum_{|u|=1} e^{-t^* V(u)}  V(u) \Big]  = \int_\r x e^{- t^* x} \mu(\dd x). $$ \end{fact}

\noindent For $x\in \r$, we write $\p_x$ for a probability distribution under which the random walk $(S_n)_{n\ge 0}$ starts at $x$. Notice that $\e[S_1]\neq 0$ under (H1) whereas $\e[S_1]=0$ under (H2).  Under (H2), $(S_n)_{n\ge 0}$ is a centered random walk with finite variance so that by Lawler and Limic \cite[Theorem 5.1.7]{lawler-limic}: 
there exists a positive constant $c$ such that for any $a, b>0$ and  $x\in (-a, b)$, \begin{equation} \label{eq:ruin} \p_x\big(T_{(-\infty, -a]} < T_{[b, \infty)}\big) \ge c \frac{b-x +1}{b+a+1} \end{equation}  

\noindent with the notation 
\begin{equation}\label{TA}
T_A:= \inf\{i\ge   1: S_i \in A\}
\end{equation}

\noindent for any Borel set $A$.\\

The  many-to-one formula can be extended to optional lines in the sense of  \cite[Section 6] {BK04}. In our applications, an optional line $\L$ will be of the form 
\begin{equation}\label{def:L}
\L = \{\xi_k,\, (\xi,k) \in \partial \T \times \mathbb{N} \hbox{ such that } \tau(\xi)=k\}
\end{equation}

\noindent where $\tau(\xi):=\inf\{k\ge 0\,:\, (V(\xi_0),V(\xi_1),\ldots,V(\xi_k)) \in A_k\}$ and for each $k\ge 0$, $A_k$ is a Borel set of $\r^{k+1}$. A consequence of the many-to-one formula is that
\begin{equation}\label{L:moment1}
\e[\#\L] = \e[e^{t^*S_{\tau}} {\bf 1}_{\{\tau<\infty\}}]
\end{equation}

\noindent where $\tau:= \inf\{k\ge 0\,:\, (S_0,S_1,\ldots,S_k) \in A_k\}$. Here and before, we used the convention that $\inf \emptyset =\infty$. Finally, consider
\begin{equation}
M_{\L} := \sum_{u\in \L} e^{-t^* V(u)}. 
\end{equation}

\noindent Its first moment is 
\begin{equation}\label{M:moment1}
\e[M_{\L}]= \p(\tau<\infty).
\end{equation}

\noindent It is a consequence of the many-to-one formula. We will need to compute its second moment. For real numbers $s_0,\ldots,s_k,s_{k+1}$, let  
\begin{equation}\label{def:p}
p(s_0,\ldots,s_k,s_{k+1}):= \p(k+1\le   \tau<\infty \mid S_0=s_0,\ldots, S_k=s_k,S_{k+1}=s_{k+1})
\end{equation}

\noindent and
\begin{equation}\label{def:psi}
\psi(s_0,\ldots,s_k):=  \e\Big[\sum_{|u|,|v|=1,u\neq v} e^{-t^*V(u)}e^{-t^* V(v)}p_up_v  \Big]
\end{equation}

\noindent where $p_u=p(s_0,\ldots,s_k,  s_k+ V(u))$,  $p_v=p(s_0,\ldots,s_k, s_k+ V(v))$. With this notation,
\begin{equation}\label{M:moment2}
\e[(M_\L)^2] = \e\left[e^{-t^*S_{\tau}} {\bf 1}_{\{\tau<\infty\}} \right]+  \e\left[\sum_{k=0}^{\tau-1} e^{-t^* S_k} \psi(S_0,\ldots,S_k)   \right].
\end{equation}

\noindent  Let us prove it. For $w\in \T$, we write $w <\L$ if $(V(w_0),V(w_1),\ldots,V(w_k)) \notin A_k$ for all $k\le |w|$. By decomposing $(M_\L)^2=\sum_{x,y \in \L} e^{-t^*V(x)}e^{-t^* V(y)}$ with respect to the most recent common ancestor $w$ of $x$ and $y$, we get
$$
(M_\L)^2
=
\sum_{w\in \L} e^{-2t^*V(w)} + \sum_{w<\L}\sum_{u\neq v \text{ children of }  w} M_{\L}^u M_{\L}^v
$$

\noindent with $M_\L^u:= \sum_{x\in \L^u} e^{-t^*V(x)}$ for any $u\in \T$, where $\L^u$ is the set of vertices $x$ descendants of $u$ which belong to $\L$. By the branching property and \eqref{M:moment1},
$$
\e\left[ \sum_{w<\L}\sum_{u\neq v \text{ children of }  w} M_{\L}^u M_{\L}^v\right] = \e\left[ \sum_{w<\L}\sum_{u\neq v \text{ children of } w} e^{-t^*V(u)}e^{-t^* V(v)} \widetilde p_u\widetilde p_v \right]
$$
 
 \noindent with $\widetilde p_u:= p(V(u_0),V(u_1),\ldots,V(u))$ and the analog for $\widetilde p_v$. Another use of the branching property shows that 
 $$
 \e\left[ \sum_{w<\L}\sum_{u\neq v \text{ children of } w } e^{-t^*V(u)}e^{-t^* V(v)}  \widetilde p_u\widetilde p_v \right]
 =
 \e\left[ \sum_{w<\L} e^{-2 t^*V(w)}  \psi(V(w_0),V(w_1),\ldots, V(w))  \right].
 $$

\noindent  With another use of the many-to-one formula,
\begin{align*}
\e\left[\sum_{w\in \L} e^{-2t^*V(w)}  \right]  & = \e\left[e^{-t^*S_{\tau}} {\bf 1}_{ \{ \tau<\infty \}} \right] , \\
 \e\left[ \sum_{w<\L} e^{-2 t^*V(w)} \psi(V(w_0),V(w_1),\ldots, V(w))  \right]  & = \e\left[    \sum_{k=0}^{\tau-1} e^{-t^* S_k}  \psi(S_0,\ldots,S_k) \right].
\end{align*}

\noindent It yields \eqref{M:moment2}.

\section{Proof of Theorem \ref{t:localtime}}
\label{s:localtime}

It is enough to prove that $\p^*$-a.s.\ for $n$ large enough, 
\begin{equation}\label{eq:Nxi}
\sup_{\xi\in \partial \T} N_{\xi}^n \le n^\kappa
\end{equation}

\noindent for any $\kappa>1$ under (H1) and $\kappa>2$ under (H2). Fix such a $\kappa$. Let $a>0$ be a small constant whose value will be determined later.  Without loss of generality, we will prove \eqref{eq:Nxi} where $N_{\xi}^n$ is rather defined as
$$
N_{\xi}^n= \sum_{\ell=0}^\infty {\bf 1}_{\{V(\xi_\ell) \in [ an,a(n+1)] \} }.
$$

\noindent For any $u\in \T$ and $k\in \z$, we let $$N_u^n:= \sum_{\ell=0}^{\vert u \vert} 1_{\{V(u_\ell) \in [a n, a(n+1)]\}}.$$ 

\noindent We first prove \eqref{eq:Nxi} under (H1). We  introduce the following set $$\L_n:=\{u\in \T: N_u^n \ge \lfloor n^\kappa\rfloor, \max_{0\le k<|u|} N_{u_k}^n < \lfloor n^\kappa\rfloor \}$$

\noindent which is the optional line of the particles in the branching random walk stopped when they first visited $\lfloor n^\kappa\rfloor $ times the set $[a n, a (n+1)]$. It corresponds to \eqref{def:L} with $\tau(\xi)$ given by
\begin{equation}\label{def:tau}
\tau_n(\xi):=\inf\{k\ge 0\,:\, N_{\xi_k}^n \ge  \lfloor n^\kappa\rfloor\}. 
\end{equation}

\noindent Note that 
\begin{equation}\label{eq:NxiLn}
\left\{\sup_{\xi\in \partial \T} N_\xi^n \ge n^\kappa \right\} \subset \{\L_n \neq \emptyset\}.
\end{equation}

\noindent  Equation \eqref{L:moment1}  yields that   \begin{equation}   \e \left(\# \L_n\right)
= \e\left(e^{t^*S_{\tau_n}} 1_{\{\tau_n<\infty\}}\right),
 \end{equation}

\noindent where 
\begin{equation}\label{def:taun}
\tau_n:= \inf\left\{i\ge 0: \sum_{\ell=0}^i 1_{\{S_\ell \in [a n, a (n+1)]\}} \ge   \lfloor n^\kappa\rfloor\right\}.
\end{equation}

\noindent We have $S_{\tau_n} \in [a n, a (n+1)]$ on $\{\tau_n<\infty\}$. If we set $\Sigma_n:= \sum_{\ell=0}^{\infty} 1_{\{S_\ell \in [a n, a(n+1)]\}}$, then 
\begin{equation}\label{eq:Ln}
\p(\L_n \neq \emptyset)
\le 
\e \left[ \# \L_n\right]
\le
e^{a t^* (n+1)} \p\Big( \Sigma_n \ge \lfloor n^\kappa\rfloor\Big).
\end{equation}

\noindent Under (H1), the random walk $(S_n)_{n\ge 0}$ either drifts to $+\infty$ a.s. or to $-\infty$ a.s. according to $\e[S_1]>0$ or $\e[S_1]<0$. Suppose for example that $S_n\to\infty$ a.s.  In particular $\p(S_n\ge 0 , \, \forall\,n\ge 0)>0$. Choose $a>0$ such that $\p(S_1>a)>0$. For any $n\ge 0$ and $x \in [a n, a (n+1) ]$, $$\p_x\Big( T_{[a n, a (n+1) ]}=\infty\Big)\ge \p(S_n> a, \,\forall\, n\ge 1) =:c>0$$ 

\noindent  with the notation \eqref{TA} for $T_{[a n, a (n+1) ]}$. By the Markov property,  for any $\ell \ge 1$, $\p(\Sigma_n \ge \ell) \le (1- c)^{\ell-1}$. In view of \eqref{eq:Ln}, using that $\kappa>1$, we get that $\sum_n \p(\L_n \neq \emptyset)$ is summable. The Borel--Cantelli lemma implies that $\L_n=\emptyset$ for $n$ large enough, and we deduce from \eqref{eq:NxiLn} that   \eqref{eq:Nxi} holds. The same reasoning deals with the case $S_n\to-\infty$.

\bigskip

It remains to prove \eqref{eq:Nxi} under (H2). We apply the same arguments using the standard strategy of introducing a barrier. Recall that under (H2), $m_n \to + \infty$ $\p^*$-a.s. Let $c>0$. Now $\L_n=\L_n(c)$ is defined as 
$$\L_n:=\{u\in \T: N_u^n \ge \lfloor n^\kappa\rfloor, \max_{0\le k<|u|} N_{u_k}^n < \lfloor n^\kappa\rfloor, \min_{0\le k\le |u|} V(u_k) > - c \}.$$

\noindent It is of the form \eqref{def:L} with $\tau(\xi)=\tau_n(\xi)$ of \eqref{def:tau} if $\min_{0\le k\le \tau_n(\xi)} V(\xi_k) > - c$ and $\tau(\xi)=\infty$ otherwise. Using \eqref{def:taun} for the definition of $\tau_n$, the many-to-one formula \eqref{L:moment1} tells us that 
$$ 
\p(\L_n \neq \emptyset) \le \e \left[\# \L_n\right]
= \e\left[e^{t^*S_{\tau_n}} 1_{\{\min_{0\le i \le \tau_n} S_i > -c\}}\right]
\le 
e^{a t^* (n+1) } \p\Big( \min_{0\le i \le \tau_n} S_i > -c \Big).
$$ 

\noindent  With the notation \eqref{TA} and   $\Sigma^{(c)}_n:= \sum_{\ell=0}^{T_{(-\infty, -c]}} 1_{\{S_\ell \in [ a n, a (n+1)]\}}$, we get 
$$
\p(\L_n \neq \emptyset)  \le e^{ a t^* (n+1)} \p\Big(\Sigma^{(c)}_n \ge \lfloor n^\kappa\rfloor\Big).
$$

\noindent We now choose  $a >0$ small enough such that $\p(S_1 < - a) >0$. By \eqref{eq:ruin}, there exists $c_1>0$ such that for any $n\ge 1$ and $x \in [a n, a (n+1)]$, $$  \p_x\Big(T_{(-\infty, -c]} < T_{[a n , a (n+1)]}\Big) \ge \frac{c_1}{n}.
$$

\noindent  This, in view of the Markov property, implies that  for any $\ell \ge 1$, $$\p(\Sigma^{(c)}_n \ge \ell) \le \left(1- \frac{c_1}{n}\right)^{\ell-1}.$$

\noindent Hence we get that $$\p(\L_n \neq \emptyset)
\le
e^{a t^*(n+1)} e^{-c_1 (\lfloor n^\kappa\rfloor -1) n^{-1}}.$$

\noindent For any $\kappa>2$, the sum over $n$ converges. The Borel-Cantelli lemma yields that almost surely, for all large $n$, $\L_n = \emptyset$, hence $\sup_{\xi\in \partial \T} N_\xi^n \le n^\kappa$ on the event  $\{\inf_{u\in\T} V(u) > -c\}$. Since $\p^*$-a.s. $\inf_{u\in\T} V(u)>-\infty$,  we conclude by letting $c\to \infty$. $\Box$

\section{Proof of Theorem \ref{t:discounted}}
\label{s:discounted}

\subsection{Proof of Part (i).} 

It is a consequence of the following lemma.

\begin{lemma}\label{l:supdiscounted}
Under the conditions of  Theorem \ref{t:discounted} (i), $\sup_{u\in \T} e^{-V(u)}\eta_u=+\infty$ \mbox{\rm $\p^*$-a.s.}
\end{lemma} 

\noindent {\it Proof of the lemma}. Let $b>0$ be small enough such that $t^*- 2b > \zeta(1+b)$. For a point process $\widehat{\aleph}\le \aleph$, let $\widehat{\mu}:=\e[\widehat{\aleph}]$,  $\widehat{\varphi}(t):=\log \int_\r e^{-t x} \widehat{\mu} (\dd x) $. Since $\varphi(t^*)=0$ and $\varphi(t)>0$ for all $t<t^*$, monotone convergence implies that we can choose $\widehat{\aleph} \le \aleph$ such that:
\begin{enumerate}[(i)]
\item there exists $a>0$ such that $\widehat{\aleph}(\r)\le a$ a.s.,
\item  the support of $\widehat{\aleph}$ is in $(-a,a)$ a.s.,
\item  $\widehat \varphi(0)>0$,
\item  $\widehat{\varphi}(s)=0$ for some $s\in (0,t^*)$ which satisfies $s - 2b > \zeta(1+b)$.  
\end{enumerate}

\bigskip

We let $(V(u),\, u\in \widehat{\T})$  be the branching random walk associated to the point process $\widehat{\aleph}$, which is contained in the original branching random walk $(V(u),\, u\in \T)$, i.e. $\widehat{\T}\subset \T$. Necessarily $\widehat \varphi(t^*)<0$. Since $s \widehat \varphi\,'(s) <0 = \widehat \varphi(s)$,  Biggins' theorem \cite{biggins-mart-cvg} implies that $\sum_{\vert u\vert =n,\, u\in \widehat{\T} } e^{-sV(u)}$ converges a.s. and in $L^1$ to some random variable $\widehat W$ which is positive on $\{\partial \widehat\T \neq \emptyset\}$.  Let 
 $$
 \widehat {\mathcal L}_k:=\{u\in \widehat \T\,:\, V(u)\ge k,\, \max_{0\le j < |u|} V(u_j)<k \}
 $$ 
 
 \noindent be the set of particles which first cross level $k$. Notice that for any $u\in \widehat {\mathcal L}_k$, $V(u)\in [k,k+a]$. Theorem 6.1  in \cite{BK04} shows that  $\sum_{u \in \widehat{\mathcal L}_k} e^{- s  V(u)}$ also converges a.s. to $\widehat W$. Since 
 $$
 \sum_{u \in \widehat{\mathcal L}_k} e^{- s V(u)}\le e^{-s k} \# \widehat {\mathcal L}_k,
 $$
 
\noindent we deduce that  a.s. on $\{\widehat W>0\}$, for $k$ large enough,
 $$
 \# \widehat {\mathcal L}_k\ge  e^{(s-b) k}.
 $$
 
 \noindent On the other hand, \eqref{taileta} implies the existence of a constant $c>0$ such that for any $x\ge 1$, $\p(\eta \le x)\le 1- cx^{-\zeta-\varepsilon}$ where $\varepsilon>0$ is such that $(1+b)(\zeta+\varepsilon)=s -2 b$. By the branching property applied to the optional line $\widehat{\mathcal L}_k$,
$$
\p\left( \forall\, u\in \widehat{\mathcal L}_k, \, \eta_u \le e^{(1+b)k},\,  \# \widehat {\mathcal L}_k\ge  e^{(s-b) k}\right) \le \left(1-c e^{-(s-2b) k}\right)^{ e^{(s-b) k} }\le e^{-c e^{bk}}.
$$

\noindent The right-hand side is summable in $k$, hence the Borel--Cantelli lemma shows that a.s. on $\{\widehat W>0\}$, for all $k$ large enough one can find $u\in \widehat {\mathcal L}_k$ such that $\eta_u > e^{(1+b)k}$. For such a $u$, $e^{-V(u)}\eta_u\ge e^{-k-a}e^{(1+b)k} = e^{bk-a}$. Therefore $\p\left(\sup_{u\in \T} e^{-V(u)}\eta_u =\infty\right)\ge \p\left(\widehat  W>0\right)>0$.  Since $\{  \sup_{u\in \T} e^{-V(u)}\eta_u =\infty\}$ is an inherited property, $\p^*(\sup_{u\in \T} e^{-V(u)}\eta_u =\infty)=1$. $\Box$

\bigskip

\subsection{Proof of Part (ii).}  Suppose now that $\zeta> t^*$. Let $r\in (t^*, \zeta)$ and $b\in (0,1)$ be such that $t^*-b r <0$. By \eqref{taileta}, there exists  a constant $c>0$ such that for all $x\ge 1$, 
 \begin{equation}\label{boundeta}
 \p(\eta \ge x)\le c x^{-r}.
 \end{equation}
 
  \noindent By  Theorem \ref{t:localtime} (already proved in Section \ref{s:localtime}), there exists $\kappa>0$ such that $\sup_{\xi} N_\xi^n \le \lfloor n^{\kappa}\rfloor$ for all $n$ large enough. Since 
 $$
D(\xi)=\sum_{n=0}^\infty e^{-V(\xi_n)} \eta_{\xi_n}\le \sum_{n=0}^\infty e^{-n} \sum_{u\in \xi\,:\, n\le V(u)<n+1} \eta_u
\le
\sum_{n=0}^\infty e^{-n} \, N_\xi^n\, \sup_{u\in \T} \eta_u {\bf 1}_{\{ n\le V(u) <n+1\}},
$$
 
 \noindent it is enough to show that for $n$ large enough, for all $u\in  \T$, 
 \begin{equation}\label{eq:ebn}
 \eta_u {\bf 1}_{\{ n\le V(u) <n+1 \}} \le e^{bn}.
 \end{equation}
    
   \noindent  For $u\in \T$ and $n\ge 0$, define  $N_u^n:= \sum_{\ell=0}^{\vert u\vert} {\bf 1}_{\{n\le V(u_\ell)< n+1\}}$. For $n\ge 0$ and $k\ge 1$, we consider the optional line
    $$
  {\mathcal L}_n(k):=\{ u\in \T \, :  \max_{0\le j<|u| }N_{u_j}^n< k,\, N_u^n =k \}.
    $$

   \noindent   It is the optional line  \eqref{def:L} associated to $\tau(\xi)=\inf\{ \ell \ge0\,:\, N_{\xi_\ell}^n =k \}$. The many-to-one formula \eqref{L:moment1} implies that  $\e[\#  {\mathcal L}_n(k)] \le e^{t^*(n+1)}$. By the branching property and the union bound, we get that 
   $$ 
\p( \exists\, u\in   {\mathcal L}_n(k)\,:\,  \eta_u\ge e^{b n})\le e^{t^*(n+1)} \p(\eta\ge e^{bn}) \le c e^{t^*(n+1) - brn} 
   $$
   
   \noindent by \eqref{boundeta}. The sum on the right-hand side over $k\le\lfloor n^{\kappa} \rfloor$ and $n\ge 0$ is finite. The Borel--Cantelli lemma implies that for $n$ large enough, $\eta_u {\bf 1}_{\{ n\le V(u) < n+1,\, N_u^n\le \lfloor n^\kappa \rfloor\}} \le e^{bn}$ for all $u\in \T$ a.s. Since $ \sup_{u\in \T} N_u^n\le \lfloor n^\kappa \rfloor$ a.s. for $n$ large, we proved \eqref{eq:ebn}. $\Box$

  \begin{remark} \label{r:continuity} Suppose that $\zeta > t^*$. For $\xi, \xi' \in \partial\T$, let $\dd_{\partial\T}(\xi, \xi'):= 2^{-|\xi \wedge \xi'|}$ be the usual ultra metric, where $|\xi \wedge \xi'|:= \sup\{n\ge 0: \xi_n= \xi'_n\}$. Recall from \eqref{def-D}  that $  D(\xi)= \sum_{n=0}^\infty e^{-V(\xi_n)}\eta_{\xi_n}$.   Since $|D(\xi)- D(\xi')|  \le \sum_{n=|\xi \wedge \xi'|+1}^\infty (e^{-V(\xi_n)}\eta_{\xi_n}+ e^{-V(\xi'_n)}\eta_{\xi'_n})$, an application of \eqref{eq:ebn}  yields that $\p^*$-a.s., $D(\xi^{(k)})\to D(\xi)$ for $k\to\infty$, if $\dd_{\partial\T}(\xi, \xi^{(k)}) \to 0$. In other words, $\p^*$-a.s., $\xi \mapsto D(\xi)$ is continuous on $\partial\T$. \end{remark}

\section{Proof of Theorem \ref{t:limsup}}
\label{s:limsup}

   Recall that $(S_n,\, n\ge 0)$ is  associated to the branching random walk through \eqref{many-to-one}. From the assumptions on $\mu$, $(S_n,\,n\ge 0)$ is integer-valued and $\p(S_1\ge 2)=0$. We introduce the local times 
   \begin{equation}\label{def:elln}
   \ell_n:=\sum_{j=0}^n {\bf 1}_{\{S_j=0\}}
   \end{equation}
   
\noindent    and the inverse local times defined by
  $\sigma_0:=0$ and for $j\ge 1$
  \begin{equation} \label{def:r}
 \sigma_j:=\inf\{k>\sigma_{j-1}: S_k=0\}
 \end{equation}

\noindent  with the convention that $\inf\emptyset=\infty$. In the case (i) (namely under (H1')),  $\e[S_1]>0$ and $\p(S_1<0)>0$. We let 
 \begin{equation}\label{def:q}
 q:=\p(\sigma_1<\infty) \in (0,1).
 \end{equation}

\noindent In the case (ii) (namely under (H2)), the random walk $(S_n,\, n\ge 0)$ is centered with finite variance and by Lemma \ref{l:hitting} there exists some $\theta>0$ such that

\begin{equation}  q(n):= \p\Big(\min_{0\le i \le \sigma_1} S_i \le - n\Big) \sim \frac{\theta}{n}, \qquad n\to\infty. \label{theta} \end{equation}

\bigskip
  
Recall that the minimum of the branching random walk goes to $\infty$,  $\p^*$-a.s.  in both cases. Notice that the sets $\{u\in \T\,:\, V(u)=k,\, \max_{0\le j<|u|} V(u_j)<k \}$ indexed by $k\ge 0$ form a Galton--Watson process.  According to \eqref{L:moment1}, its  mean offspring is $e^{t^*}$. In particular, the line  $\{u\in \T\,:\, V(u)=n,\, \max_{0\le j<|u|} V(u_j)<n \}$ has size $e^{t^* n +o(n)}$ as $n\to \infty$ $\p^*$-a.s. 
   By considering the descendants of these particles, we see that in order to prove that  $\sup_{\xi\in \partial \T} N_\xi^n \le d_n$ as $n\to \infty$, where     $$ d_n:= \begin{cases}  d \, n , \qquad &\mbox{in the case {\rm (H1')}}, \\
d\, n^2,   \qquad & \mbox{in the case {\rm (H2)}}, \end{cases}  \qquad \mbox{ with } d>0, $$

\noindent it is enough to show that  there exists $\varepsilon_d >0$ such that for $n$ large enough
\begin{equation}\label{upper_limsup}
\p(N_* >  d_n) \le e^{-t^*n -\varepsilon_d  n}
\end{equation}

\noindent where $N_*:=\sup_{\xi\in \partial T} N_\xi^0$.  Indeed, it would be a consequence of the Borel--Cantelli lemma. Similarly, to prove that $\sup_{\xi\in \partial \T} N_\xi^n \ge d_n$ as $n\to \infty$, it is enough to show that  there exists $\varepsilon'_d>0$ such that for $n$ large enough
\begin{equation}\label{lower_limsup}
\p(N_* > d_n) \ge e^{-t^*n +\varepsilon'_d n}.
\end{equation}

The proof  of Theorem \ref{t:limsup} is presented  in the following three subsections, separately for the case (H1') and  (H2).   For $u\in \T$ and $k\ge 0$, we let $N_u^k:=\sum_{i=0}^{|u|} 1_{\{V(u_i)=k\}}$.

   \subsection{Proof of Theorem \ref{t:limsup} (i)}

  \noindent  Let $d  > -\frac{t^*}{\log q}$. We want to prove \eqref{upper_limsup}. For notational brevity, we treat $d n$ as an integer. We consider the optional line 
 $$ 
    \L_n:=\{u\in\T:  N_u^0 > d  n , \max_{0\le j<|u|}  N_{u_j}^0 \le d  n\}
  $$
  
  \noindent which is the set of particles stopped when their local time at level $0$ exceeds $d  n$. It is of the form \eqref{def:L} with $\tau(\xi)$ given by
  \begin{equation}\label{def:tau_upper_H1}
  \sigma_{d  n}(\xi):=\inf\{k\ge 0\,:\, N^0_{\xi_k}> d  n\}.
  \end{equation}
  
\noindent   The many-to-one formula \eqref{L:moment1} implies that
  $$
  \e[\#\L_n] = \p\left(\sigma_{d  n}<\infty \right) = q^{d  n}
  $$
  
  \noindent  with the notation \eqref{def:r}. Then \eqref{upper_limsup} follows from $\p(N_* > d  n)\le \p(\L_n\neq \emptyset)\le \e[\#\L_n]$.

  We turn to the lower bound. Let $d  \in (0,-\frac{t^*}{\log q})$. We want to prove \eqref{lower_limsup}. Take $K\ge 1$ large enough such that $e^{t^*}q_K^d >1$ where
  $$
  q_K:=\p(\sigma_1<\infty,\, \min_{0\le n\le  \sigma_1} S_n \ge -K).
  $$
  
  \noindent We let 
$$
\L_n(K):= \{u\in\T:  N_u^0 > d n , \max_{0\le j<|u|}  N_{u_j}^0 \le d  n, \,\min_{0\le j\le |u|}  V(u_j)\ge -K\}
$$  
  
 \noindent so that if $\L_n(K)$ is nonempty then $N_*> d  n$. Using the notation \eqref{def:tau_upper_H1}, $\L_n(K)$ is of the form \eqref{def:L} with $\tau(\xi)= \sigma_{d  n}(\xi)$ if $ \min_{0\le j \le \sigma_{d  n}(\xi)}  V(\xi_j) \ge -K$, and $\tau(\xi)=\infty$ otherwise. We introduce 
 $$
 M_{\L_n(K)}:= \sum_{u\in \L_n(K)} e^{-t^* V(u)}.
 $$ 
  
  \noindent We have by \eqref{M:moment1} 
  \begin{equation}\label{eq:moment1_upper_H1}
  \e[M_{\L_n(K)}] = \p(\tau<\infty)=\p(\sigma_{d  n}<\infty,\, \min_{0\le j\le \sigma_{d  n}} S_i\ge -K)=q_K^{d  n}
  \end{equation}
  
\noindent   where now $\tau= \sigma_{d  n}$ if $ \min_{0\le j \le \sigma_{d  n}}  S_j \ge -K$, and $\tau=\infty$ otherwise. Recall the notation $p(s_0,\ldots,s_{k+1})$ and $\psi(s_0,\ldots,s_k)$ in \eqref{def:p} and \eqref{def:psi}, which read in our setting 
\begin{align*}
p(s_0,\ldots,s_k,s_{k+1})  & =\p(\tau<\infty \mid S_0=s_0,\ldots, S_k=s_k,S_{k+1}=s_{k+1}) \\
& = 
\p(  \sigma_{d  n}<\infty,\, \min_{0\le j\le \sigma_{d  n}} S_i\ge -K \mid S_0=s_0,\ldots, S_k=s_k,S_{k+1}=s_{k+1}) 
\end{align*}

\noindent and
$$
\psi(s_0,\ldots,s_k)=  \e\Big[\sum_{|u|,|v|=1,u\neq v} e^{-t^*V(u)}e^{-t^* V(v)}p_u p_v  \Big]
$$

\noindent where $p_u=p(s_0,\ldots,s_k, s_k+ V(u))$,  $p_v=p(s_0,\ldots,s_k,  s_k+ V(v))$.  The second moment of $M_{\L_n(K)}$ is given by \eqref{M:moment2}, i.e.
\begin{equation}\label{M:moment2_proofH1}
\e[(M_{\L_n(K)})^2] = \e\left[e^{-t^*S_{\tau}} {\bf 1}_{\{\tau<\infty\}} \right]+  \e\left[\sum_{k=0}^{\tau-1} e^{-t^* S_k} \psi(S_0,\ldots,S_k)   \right].
\end{equation}

\noindent Notice that   $$
  \psi(S_0,\ldots,S_k){\bf 1}_{\{k<\tau\}}=\psi(S_0,\ldots,S_k){\bf 1}_{\{k<\sigma_{d  n},\, \min_{0\le j\le k} S_j \ge -K\}}. 
  $$

\noindent   Indeed,   $p(s_0,\ldots,s_{k+1})=\p(k+1\le \tau<\infty \,|\, S_0=s_0,\ldots,S_{k+1}=s_{k+1})={\bf 1}_{\{\sigma_{d n} \ge k+1\}} \p(\sigma_{d n}<\infty,  \min_{0\le j\le \sigma_{d n}} S_j \ge -K \,|\, S_0=s_0,\ldots,S_{k+1}=s_{k+1})$. Therefore $\psi(S_0, ..., S_k)=0$ on $\{ \sigma_{d n} \le k\} \cup\{\min_{0\le j\le k} S_j < -K\}$, and  the above equality follows.

In \eqref{M:moment2_proofH1},
$$
\e\left[e^{-t^*S_{\tau}} {\bf 1}_{\{\tau<\infty\}} \right] = \p(\tau<\infty)=q_K^{d  n}
$$

\noindent by \eqref{eq:moment1_upper_H1}. The last term of  of  \eqref{M:moment2_proofH1} is
  \begin{equation}\label{upper_H1_proof}
  \e\left[\sum_{k=0}^{\sigma_{d  n} -1} e^{-t^* S_k} \psi(S_0,\ldots,S_k) {\bf 1}_{\{\min_{0\le j\le k} S_j \ge -K\}} \right].
  \end{equation}
  
  \noindent We bound $p(s_0,\ldots,s_k,s_{k+1})$ for $s_0,\ldots,s_k$ such that  $\ell:=\sum_{j=0}^k {\bf 1}_{\{s_j=0\}} \le d n$ and $\min_{0\le j\le k} s_i\ge -K$. By the strong Markov property  at the first hitting time of $0$ after time $k + 1$,  we see that $p(s_0,\ldots,s_k,s_{k+1}) \le q_K^{d  n-\ell}$. Hence $\psi(s_0,\ldots,s_k) \le c'q_K^{2(d n-\ell)}$ with $$c'=\e\Big[ \sum_{|u|,|v|=1,u\neq v} e^{- t^* V(u)}e^{- t^* V(v)} \Big].$$  Therefore the expectation in \eqref{upper_H1_proof} is less than
  $$
  c' q_K^{2d n}  \e\left[\sum_{k=0}^{\sigma_{d n} -1} e^{- t^* S_k }  q_K^{-2\ell_k}{\bf 1}_{\{\min_{0\le j\le k} S_j \ge -K\}} \right]
  $$ 
  
\noindent with the notation \eqref{def:elln}.  Discussing over the value of $\ell_k$ and using the strong Markov property at $\sigma_{\ell-1}$ below, we obtain that 
\begin{align*} 
\e\left[\sum_{k=0}^{\sigma_{d n} -1}  e^{- t^* S_k } q_K^{-2 \ell_k}{\bf 1}_{\{\min_{0\le j\le k} S_j \ge -K\}} \right] 
&=   \sum_{\ell=1}^{d n}  q_K^{-2 \ell} \e \left[\sum_{k=\sigma_{\ell-1}}^{\sigma_\ell-1} e^{- t^* S_k } {\bf 1}_{\{\min_{0\le j\le k} S_j \ge -K, \sigma_\ell< \infty\}} \right] 
\\
&=   c''\, \sum_{\ell=1}^{d n}  q_K^{-2 \ell} \p\left(\sigma_{\ell-1}< \infty,  \min_{0\le j\le \sigma_{\ell-1}} S_j \ge -K \right)
\\
&= 
 c'' \, \sum_{\ell=1}^{d n}  q_K^{-  (\ell+1)},\end{align*}

\noindent with $$c'':= \e \left[\sum_{k=0}^{\sigma_1-1} e^{- t^* S_k } {\bf 1}_{\{\min_{0\le j\le k} S_j \ge -K, \sigma_1< \infty\}} \right] .$$

\noindent By Lemma \ref{l:renewal1}, $c''<\infty$. We get that the second moment of $M_{\L_n(K)}$ is bounded by $c''' q_K^{d n} $. Recall  \eqref{eq:moment1_upper_H1}. Since
$$
\p(N^* > d n)\ge \p(\#\L_n(K) >0) = \p(M_{\L_n}(K) >0)\ge \frac{\e\left[M_{\L_n(K)} \right]^2}{\e\left[(M_{\L_n(K)})^2\right]} \ge  \frac1{c'''} q_K^{d n},
$$
  
  \noindent and $q_K^d >e^{-t^*}$, we proved \eqref{lower_limsup} and then the lower bound in Theorem \ref{t:limsup} (i). $\Box$

  \subsection{Proof of Theorem \ref{t:limsup} (ii). Upper bound.}
  
 Let  $d > \frac{t^*}{2\theta}$. We show \eqref{upper_limsup}. We consider a family of continuous non-increasing   functions $f_n: \r_+\to \r_+$ such that $f_n(0) = n$, $f_n(t) >0$ for $0\le t < d n^2$ and  $f_n(t)=0$ for all $t\ge d n^2$. 
We introduce the optional line  \begin{equation}   \L_n:=\{u\in\T:  V(u) \le  - f_n(N_u^0), \forall\, 0\le j<|u|, \,  V(u_j) > - f_n(N_{u_j}^0)\}  \label{Ln} \end{equation}
\noindent where we have stopped the particles $u$ which go below $- f_n(N_u^0)$ for the first time.  It is the optional line \eqref{def:L} associated to $\tau(\xi)$ which is 
\begin{equation}\label{def:taun_upper_H2}
\tau_n(\xi):=\inf\{k\ge 0\,:\, V(\xi_k) \le -f_n(N_{\xi_k}^0)\}.
\end{equation}
 We notice that $\{N_*> d n^2\} \subset \{\L_n\neq \emptyset\}$  as $f_n(d n^2)=0$ so that $ \p\left( N_* > d n^2\right) 
 \le
 \e [\# \L_n]$.  It follows from \eqref{L:moment1} that  \begin{equation}    \p\Big( N_* > d n^2\Big) \le
\e \Big[ e^{t^*S_{\tau_n}}  \Big]
=:b_n,  \label{b_n}\end{equation}

\noindent  where   $\tau_n:= \inf\{i\ge 0: S_i \le - f_n(\ell_i)\}$ with the notation \eqref{def:elln} for the local time process $(\ell_i)$. Observe that $\tau_n \le \sigma_{\lceil d n^2\rceil}$, with $ \sigma_j $ defined in \eqref{def:r}.  In the case (ii), $(S_n)$ is a centered random walk with finite variance, so $\tau_n$ is finite almost surely.

To estimate $b_n$, we discuss on the value of $\ell_{\tau_n}$. Note that 
\begin{equation}\label{eq:bn}
 b_n \le \e \left[ e^{-  t^*f_n(\ell_{\tau_n})} \right] = \sum_{k=1}^{\lceil d n^2\rceil} e^{- t^*f_n(k)} \p(\ell_{\tau_n}=k)
\end{equation}
\noindent 
where we used the fact that  $\ell_{\tau_n} \le {\lceil d n^2\rceil}$.   Recall the notation $q(n)$ in \eqref{theta}. We observe that, for any $k\ge 1$,
\begin{align}
\p(\ell_{\tau_n}=k) &= \p \Big( \forall 1\le j < k, \min_{\sigma_{j-1}\le i \le \sigma_{j}} S_i > - f_n(j),  \min_{\sigma_{k-1}\le i \le \sigma_k} S_i \le - f_n(k)\Big)\nonumber\\
&=
 \prod_{j=1}^{k-1} (1-q(f_n(j))) \, q(f_n(k))\label{prodq}
\end{align}

\noindent with the convention  $\prod_{\emptyset}=1$. Hence, for all $k\ge 1$,
\begin{equation}\label{pnell}
e^{- t^* f_n(k) }\p(\ell_{\tau_n}=k) \le  e^{- t^*f_n(k) - \sum_{j=1}^{k-1} q(f_n(j))}.
\end{equation}

\noindent  As $d > \frac{t^*}{2\theta}$, we may choose and fix some $\theta'\in (0, \theta)$ such that $\alpha:= \sqrt{2 d  \frac{\theta'}{t^*}} >  1$. By \eqref{theta}, there exists $k_0\ge 1$ such that for all $k\ge k_0$, $q(k)\ge \frac{\theta'}{k}$. Let $f(t):= \alpha \sqrt{(1-\frac{t}{d})_+}$ and $f_n(t):= n f(t/n^2)$ for $t\ge 0$. We check that for any $t\le d n^2$,
\begin{equation}\label{integral}
t^*f_n(t) + \theta' \int_0^{t} \frac{\dd s}{f_n(s)}= \alpha t^*n.
\end{equation}

\noindent It implies that $ t^*f_n(k) + \theta' \sum_{j=1}^{k+1} \frac1{f_n(j)} \ge \alpha t^* n$ for any $k\ge 0$. Let $m_n$ be the minimal integer such that $f_n(m_n+2)<k_0$. We can check that $f_n(m_n)$ is bounded in $n$. For $k \le m_n$, we get
$$
t^* f_n(k) + \sum_{j=1}^{k-1} q(f_n(j)) \ge t^* f_n(k) + \theta' \sum_{j=1}^{k-1} \frac{1}{f_n(j)} \ge   \alpha t^* n - \frac{2\theta'}{k_0}
$$

\noindent while for $k> m_n$, we write
$$
t^* f_n(k) + \sum_{j=1}^{k-1} q(f_n(j)) \ge \sum_{j=1}^{m_n-1} q(f_n(j)) \ge \alpha t^* n   - \frac{2\theta'}{k_0}-t^*f_n(m_n).
$$

\noindent In view of \eqref{pnell}, equation \eqref{eq:bn} becomes $b_n \le e^{-\alpha t^* n}\lceil d n^2 \rceil e^{  \frac{2\theta'}{k_0}+t^*f_n(m_n)}$, proving \eqref{upper_limsup}. $\Box$

  \subsection{Proof of Theorem \ref{t:limsup} (ii). Lower bound.}

 Let $0< d <   \frac{t^*}{2\theta}$. For $\varrho \in (d,   \frac{t^*}{2\theta})$,  let $\alpha:=\sqrt{2 \varrho  \frac{\theta}{t^*}}$ and take $f_n(t)=nf(t/n^2)$ with $f(t)=\alpha\sqrt{(1-\frac{t}{\varrho})_+}$ for $t\ge 0$. Note that $\alpha <1$.
 We set $d_n:=\lceil d n^2\rceil+1$ for concision. Notice from \eqref{theta} and \eqref{integral} that uniformly in $k \le d_n$,
\begin{equation}\label{est_prodq}
e^{-t^*f_n(k)}\prod_{j=0}^{ k-1} (1-q(f_n(j)))=e^{-\alpha t^* n + o(n)}.
\end{equation}

\noindent We consider a slightly modified version of  the optional line $\L_n$ defined in \eqref{Ln}: 
   Let $(\lambda_n)_{n\ge0}$ be a sequence of positive real numbers such that $\lambda_n=o(n)$ as $n\to\infty$.   
   Consider 
\begin{align*}     {\L_n^{(d)}}:=\{u\in\T: & V(u) \le  - f_n(N_u^0), |V(u) + f_n(N_u^0)|\le \lambda_n, \,  
\\ &  \forall \, 0\le j<|u|,\,   V(u_j) > - f_n(N_{u_j}^0), N_u^0= d_n \} .  \end{align*}

\noindent In other words, we only select those $u\in \L_n$ such that $N_u^0 = d_n $ and the overshoot $|V(u) + f_n(N_u^0)|$ is less than $\lambda_n$. If ${\L_n^{(d)}}$ is not empty, then $N_* \ge d_n > d n^2$. We use the definition of $\tau_n(\xi)$ in \eqref{def:taun_upper_H2}. The optional line $\L_n^{(d)}$ is of the form \eqref{def:L} with $\tau(\xi)$ given by $\tau(\xi):=\tau_n(\xi)$ if $N_{\xi_{\tau_n(\xi)}}^0=d_n$  and $|V(\xi_{\tau_n(\xi)})+f_n(d_n)| \le \lambda_n$, and $\tau(\xi):=\infty$ otherwise. Let $$M_{\L_n^{(d)}}:= \sum_{u \in {\L_n^{(d)}}} e^{- t^* V(u)}.$$

\noindent   We will use \eqref{M:moment1}, with the notation  $\tau=\tau_n$ on the event $\{\ell_{\tau_n}= d_n, |S_{\tau_n}+f_n(d_n)|\le \lambda_n\}$,  and $\tau=\infty$ otherwise, where as before,    $\tau_n:= \inf\{i\ge 0: S_i \le - f_n(\ell_i)\}$.   We then get that   $ \e\left[M_{\L_n^{(d)}}\right]
=  \p \Big( \ell_{\tau_n} = d_n, |S_{\tau_n}+f_n(d_n)|\le \lambda_n\Big)$, and   
 by the strong Markov property at $\sigma_{d_n}$, 
$$
\p(\ell_{\tau_n}= d_n, |S_{\tau_n}+f_n(d_n)|\le \lambda_n  ) = \prod_{j=0}^{ d_n -1} (1-q(f_n(j))) \, x_n,
$$

\noindent  where $$x_n:=\p\left(T'_{f_n(d_n)}< \sigma_1, \, |S_{T'_{f_n(d_n)}} +  f_n(d_n)|\le \lambda_n \right)  $$  

\noindent with $T'_{f_n(d_n)}:= \inf\{i\ge 0: S_i \le - f_n(d_n)\}$. By \eqref{eq:thetaKn}, if we choose $\lambda_n:= K_{\lfloor f_n(d_n)\rfloor}$ (which we will), then $\lambda_n=o(n)$ and $x_n\sim \frac{\theta}{f_n(d_n)} \sim q(f_n(d_n))$ as $n\to\infty$. Then by \eqref{prodq} and \eqref{est_prodq}, we get that  \begin{equation}\label{elltaun}\e\left[M_{\L_n^{(d)}}\right] = e^{-\alpha t^* n +t^* f_n(d_n) + o(n)}.
 \end{equation}

   We compute now the second moment of $M_{\L_n^{(d)}}$. Recall \eqref{M:moment2} and the notation $p(s_0,\ldots,s_{k+1})$ and $\psi(s_0,\ldots,s_k)$ in \eqref{def:p} and \eqref{def:psi}.   Note that $p(s_0,\ldots,s_{k+1})=\p(k+1\le \tau<\infty \,|\, S_0=s_0,\ldots,S_{k+1}=s_{k+1})={\bf 1}_{\{\tau_n \ge k+1\}} \p( \ell_{\tau_n}= d_n , \, |S_{\tau_n}+f_n(d_n)|\le \lambda_n\, |\, S_0=s_0,\ldots,S_{k+1}=s_{k+1})$, one has $\psi(S_0, ..., S_k)=0$ on $\{ \tau_n\le k\} \cup\{\ell_k> d_n\}$. Therefore

$$
 \e\left[\sum_{k=0}^{\tau-1} e^{- t^* S_k} \psi(S_0,\ldots,S_k)   \right]
 = \e\left[\sum_{k=0}^{\tau_n-1} {\bf 1}_{\{\ell_k \le d_n\}}e^{- t^* S_k} \psi(S_0,\ldots,S_k)   \right].
$$

\noindent Let us compute $p(s_0,\ldots,s_{k+1})$. Let $\ell=\sum_{j=0}^k {\bf 1}_{\{s_j=0\}}$.  Suppose that $\ell\le  d_n -1$. By the strong Markov property of $S$ at the first hitting time of $0$ after time $k+1$, $p(s_0,\ldots,s_{k+1})$ is less than
$$
  \p \Big( \forall 1\le j< d_n(\ell), \min_{\sigma_{j-1}\le i \le \sigma_j} S_i > - f_n(j+\ell),  \min_{\sigma_{d_n(\ell)-1} \le i \le \sigma_{d_n(\ell)}} S_i \le - f_n(d_n)\Big)
  $$
  
  \noindent with $d_n(\ell):=d_n-\ell$. It is $\prod_{j=\ell+1}^{d_n-1} (1-q(f_n(j))) \, q(f_n( d_n))$.  We deduce that 
$$
\psi(s_0,\ldots,s_k) \le c' \left(\prod_{j=\ell+1}^{d_n-1} (1-q(f_n(j))) \, q(f_n( d_n))  \right)^2
$$

\noindent with $c':=\e[ \sum_{|u|,|v|=1,u\neq v} e^{- t^* V(u)}e^{- t^* V(v)} ] $. Therefore, in view of \eqref{est_prodq},  uniformly over $s_0,\ldots,s_k$,
$$
\psi(s_0,\ldots,s_k)\le   e^{2 t^* (f_n(d_n)-f_n(\ell)) +o(n)}.
$$

\noindent  When $\ell=d_n$, we use the bound $\psi(s_0,\ldots,s_k) \le c'$ since $p(s_0,\ldots,s_{k+1})$ is always smaller than $1$. We get
$$
  \e\left[\sum_{k=0}^{\tau-1} e^{- t^* S_k} \psi(S_0,\ldots,S_k)   \right]\\
  \le
 e^{o(n)}\,  \e\left[\sum_{k=0}^{\tau_n-1} {\bf 1}_{\{ \ell_{k}\le d_n  \}}e^{- t^* S_k} e^{2t^*(f_n(d_n)-f_n(\ell_k))} \right]  .
$$

\noindent We discuss on the local times $\ell_k$ to obtain
$$
 \e\left[\sum_{k=0}^{\tau_n-1} {\bf 1}_{\{ \ell_{k}\le d_n  \}} e^{- t^* S_k} e^{2t^* (f_n(d_n)-f_n(\ell_k))} \right]
 =
 \sum_{\ell=1}^{d_n} \e\left[  \sum_{k=\sigma_{\ell-1}}^{\sigma_{\ell}} {\bf 1}_{\{k<\tau_n\}} e^{- t^* S_k}  \right] e^{2 t^* (f_n(d_n)-f_n(\ell)) }.
$$

\noindent By the strong Markov property,
\begin{align*}
 \e\left[  \sum_{k=\sigma_{\ell-1}}^{\sigma_{\ell}} {\bf 1}_{\{k<\tau_n\}} e^{- t^* S_k}  \right] 
 &\le \p({\tau_n}\ge \sigma_{\ell-1}) \e\left[  \sum_{k=0}^\infty {\bf 1}_{\{\min_{0\le i\le k} S_i \ge -f_n(\ell)\}} e^{- t^* S_k}  \right] 
 \\
 &\le 
 c'' \p({\tau_n}\ge \sigma_{\ell-1})e^{t^* f_n(\ell)} 
\end{align*}

\noindent with $c'':=\sup_{x\ge 0} \e\left[  \sum_{k=0}^\infty {\bf 1}_{\{\min_{0\le i\le k}S_i +x \ge 0\}} e^{- t^* S_k- t^* x}  \right] $ which is finite by  Lemma \ref{l:renewal1}. Since $ \p(\tau_n\ge \sigma_{\ell-1})=  \prod_{j=1}^{\ell-1} (1-q(f_n(j)))$, we obtain by \eqref{est_prodq}
$\p({\tau_n}\ge \sigma_{\ell-1})e^{t^* f_n(\ell)} = e^{-\alpha t^* n + 2 t^* f_n(\ell)  + o(n)}$ hence 
$$
\e\left[\sum_{k=0}^{\tau-1} e^{- t^* S_k} \psi(S_0,\ldots,S_k)   \right]
\le
   e^{-\alpha t^* n+ 2 t^*f_n(d_n) +o(n)}.
$$

\noindent Moreover   $\e\left[e^{-t^*S_{\tau}} {\bf 1}_{\{\tau<\infty\}} \right]=\e\left[e^{-t^*S_{\tau_n}} {\bf 1}_{\{\ell_{\tau_n}=d_n, \, |S_{\tau_n}+f_n(d_n)|\le \lambda_n\}} \right] \le
e^{t^* f_n(d_n)+ \lambda_n} \p(\ell_{\tau_n}=d_n)$. By \eqref{prodq} and \eqref{est_prodq}, $\p(\ell_{\tau_n}=d_n)= e^{-\alpha t^* n +t^* f_n(d_n) + o(n)}. $  Recall that $\lambda_n=o(n)$. 
  Therefore \eqref{M:moment2} implies $\e[( M_{\L_n^{(d)}})^2] \le  e^{-\alpha t^* n + 2 t^* f_n(d_n) +o(n)}$. From 
$$
   \p(N_* > d n^2) \ge \p(M_{\L_n^{(d)}} >0) \ge \frac{(\e [M_{\L_n^{(d)}}])^2}{\e[M^2_{\L_n^{(d)}}]},
$$

\noindent we get $\p(N_* > d n^2) \ge e^{-\alpha t^*n +o(n)}$ which yields \eqref{lower_limsup} since $\alpha<1$. $\Box$

\bigskip
 
We deduce from \eqref{upper_limsup} and  \eqref{lower_limsup} the following tail behavior  of $N_*=\sup_{\xi\in \partial\T} N_\xi^0$:  
\begin{corollary} Under the assumptions of Theorem \ref{t:limsup}, we have 

\noindent
(i) in the case {\bf (H1')}, $$\lim_{n\to \infty}\frac1{n} \log \p(N_* > n) = \log q,  $$
\noindent where $q\in (0,1)$ is given by \eqref{def:q}. \\
\noindent
(ii) in the case {\bf (H2)},    $$\lim_{n\to \infty}\frac1{n^{1/2}}\log \p(N_* > n) =   - \sqrt{2 \theta t^*}, \qquad \mbox{\rm $\p^*$-a.s.}, $$ 
\noindent where $\theta$ is given by \eqref{theta}.
\end{corollary}

  \section{Proof of Theorem \ref{t:consistent}}\label{s:consistent}
  
 Recall that   $\p^*$-a.s., $\min_{|u|=n} V(u) \to \infty$. 
Let $\Theta:= \inf_{\xi\in \partial \T} \limsup_{n\to\infty} \frac{V(\xi_n)}{\sqrt{n\log n}} .$ For any $a>0$, the event $\{\Theta \ge   a \}$ is inherited in the sense that $\{\Theta \ge a\}= \cap_{|u|=1} \{\Theta^{(u)}\ge  a\}$, where $\Theta^{(u)}$ is defined exactly as $\Theta$ but for the branching random walk indexed by the subtree $\T$ rooted at $u$. Hence (see \cite{LPbook}) $\p^*(\Theta\ge a)\in \{0, 1\}$. It follows that $\p^*$-a.s., $\Theta$ is a constant. We only need to check that $\Theta \in (0, \infty)$, which will be done separately in the next two subsections. 
  
  \subsection{Lower bound: Proof of $\Theta>0$.} 
  
  Let $a>0$ whose value  will be determined later. Fix $K>0$. We will show that for any $n_0\ge 1$, \begin{equation} \label{low:Thm1.2} \lim_{n\to\infty} \p\Big( \exists |u|=n: -K \le V(u_\ell) \le a \sqrt{\ell \log \ell} , \forall\,  n_0 \le \ell \le n\Big) =0.
  \end{equation}

In fact, on the event $\{\inf_{u\in \T} V(u) \ge -K\}$, \eqref{low:Thm1.2} yields that a.s. for all $n_0$, there is no $\xi\in \partial\T$ such that $V(\xi_\ell) \le a \sqrt{\ell \log \ell} , \forall\,  n_0 \le \ell $, which a fortiori yields that $\Theta\ge a$. Since  $\p^*$-a.s.  $\inf_{u\in \T} V(u) > -\infty$, we get $\Theta\ge a$.

  It remains to show \eqref{low:Thm1.2}.   Let $C>0$, $ L>1$ and   $ r_i:= n_0+ \lfloor C i^2 \log i\rfloor$ for $1\le i \le L$, and  consider  the set  $$\L_{L, (r_i)}:= \{|u|=r_L:  -K \le V(u_{r_i}) \le a \sqrt{r_i \log r_i} , \forall  \, 1\le i \le L\}.$$
  
  Note that the probability term in \eqref{low:Thm1.2} is less that $\p(\#\L_{L, (r_i)}>0) $ for any $L$ (indeed we will let $L\to\infty$). By the Markov inequality and the many-to-one formula \eqref{many-to-one}, \begin{align*}  \p(\#\L_{L, (r_i)}>0) 
  &\le  \e[\#\L_{L, (r_i)}]
  \\
  &= 
  \e \Big[ e^{t^* S_{r_L}} \, {\bf 1}_{\{ -K \le S_{r_i} \le a \sqrt{r_i \log r_i} , \,  \forall 1\le i \le L\}}\Big] 
  \\
  &\le 
  e^{a  t^*\sqrt{r_L \log r_L}}\,  \p\Big(\cap_{1\le i \le L}  \{-K \le S_{r_i} \le a \sqrt{r_i \log r_i}\}\Big).
  \end{align*}

  \noindent By the Markov property for the random walk $S$, 
 \begin{equation} \p\Big(\cap_{1\le i \le L}  \{-K \le S_{r_i} \le a \sqrt{r_i \log r_i}\}\Big) 
  \le
  \prod_{i=1}^L \p \Big( S_{r_{i+1}-r_i} \le a \sqrt{r_i \log r_i}+K\Big). \label{Sri}\end{equation}

Let us recall some useful estimates on the random walk $S$. Under (H1'), $\e(S_1)= \int_\r x e^{- t^* x} \mu (\dd x)=: {\tt m}>0$. By assumption, $c_{\alpha_1}:=\e(|S_1|^{\alpha_1} )= \e\left[\sum_{|u|=1} e^{-t^*V(u)} |V(u)|^{\alpha_1} \right]  < \infty$ for some $\alpha_1>1$. Applying \cite[Lemma 2.1]{DDS08}, we get some constant $c'>0$ such that for all $n\ge 1$, $y\ge n^{\max(1/\alpha_1, 1/2)}$ and $x>0$, $$ \p\Big(S_n \le {\tt m} n -x , \min_{1\le i \le n} (S_i-S_{i-1}) \ge -y\Big) \le c' e^{-x/y}.$$

Using $ \p(S_1 \le -y)  \le c_{\alpha_1} y^{-\alpha_1},$ for any $y>0$, we get that  \begin{equation} \p\Big(S_n \le {\tt m} n -x\Big) \le   c' e^{-x/y} + c_{\alpha_1} \, n \, y^{-\alpha_1}. \label{rw-H1}\end{equation}

  We apply \eqref{rw-H1} to the probability terms on the right-hand side of \eqref{Sri}. Since $r_{i+1}- r_i \sim 2 C i \log i$ and $a \sqrt{r_i \log r_i} \sim a \sqrt{C} i \log i$, we will choose $C>0$ such that $$C {\tt m} > a \sqrt{C}.$$
  
   For all large $i$,   $a \sqrt{r_i \log r_i} +K\le \frac{\tt m}{2}  (r_{i+1}- r_i)$. Applying \eqref{rw-H1} to $n= r_{i+1}-r_i, x= \frac{\tt m}{2}  (r_{i+1}- r_i)$ and $y = (r_{i+1}-r_i)^{1-\delta}$ with $0< \delta < \min(\frac12, 1-\frac1{\alpha_1}), $ we have for all large $i$, say $i\ge i_1$,  $$
  \p\Big( S_{r_{i+1}-r_i} \le a \sqrt{r_i \log r_i}+K\Big)
  \le
  \p\Big( S_{r_{i+1}-r_i} \le \frac{\tt m}{2}  (r_{i+1}- r_i)\Big)
  \le
    (r_{i+1}-r_i)^{- \delta'},$$
  
  \noindent where $\delta'\in (0,  \alpha_1(1-\delta) -1 )$ is a positive constant. It follows that $$\p(\#\L_{L, (r_i)}>0) 
  \le e^{a  t^*\sqrt{r_L \log r_L}} e^{- \delta' \sum_{i=i_1}^L \log (r_{i+1}-r_i) }.$$

  \noindent Note that $\sum_{i=i_1}^L \log (r_{i+1}-r_i) \sim L\log L$ and $a  t^*\sqrt{r_L \log r_L} \sim a t^* \sqrt{C} L \log L$. Now we choose an arbitrary constant $a >0$ such that $a^2 t^*< \delta' {\tt m}$, and we can find $C>0$ such that  $C {\tt m} > a \sqrt{C}$ and $a t^* \sqrt{C} < \delta'$, so that $ e^{a  t^*\sqrt{r_L \log r_L}} e^{- \delta' \sum_{i=i_1}^L \log (r_{i+1}-r_i)} \to 0$ as $L\to \infty$. This proves \eqref{low:Thm1.2}. $\Box$
  
  \subsection{Upper bound: Proof of $\Theta< \infty$.} We want to construct an infinite ray $\xi\in \partial\T$ such that $\limsup_{n\to\infty} \frac{V(\xi_n)}{\sqrt{n \log n}} < \infty$. The idea is to construct a (inhomogeneous) Galton-Watson process $( \# {\mathcal Z}_i)_{i\ge 0}$ which will survive with positive probability. Let $r_1> 1$ be a large integer and consider $r_i:= r_1+ \lfloor i^2 \log i\rfloor$ for $i\ge 2$. Let ${\mathcal Z}_1=\{u\in \T: |u|=r_1\}$ be the set of particles alive in the $r_1$-th generation of $\T$. Let $C>0$ whose value will be determined later. For $i\ge 2$, we construct the set ${\mathcal Z}_i$ recursively by letting \begin{align*} {\mathcal Z}_i&:= \cup_{v\in {\mathcal Z}_{i-1}}   \Big\{u\in \T_v, |u|=r_i:    V(u_k) - V(v) \le C (k-r_{i-1}) , \forall r_{i-1} \le k \le r_i, 
  \\
  &  \qquad V(u_{r_i-1}) -V(v) \ge 2 C \log i,\,  V(u)-V(v)  \le C \log i \Big\},
  \end{align*}
  
 \noindent where $\T_v$ denotes the set of descendants of $v$.
  We note that for each $i$,   ${\mathcal Z}_i$ is a subset of $\{u\in \T: |u|=r_i\}$. By the branching property, $({\mathcal Z}_i)$ constitues an inhomogeneous  Galton-Watson tree. If $({\mathcal Z}_i)$ survives with positive probability,   on the survival set of $({\mathcal Z}_i)$,  we may find an infinite ray for the system $({\mathcal Z}_i)$, say $\xi$. Then $\xi$ can be also considered as an infinite ray of $\partial\T$. Moreover  by definition of $({\mathcal Z}_i)$, $V(\xi_{r_i}) -V(\xi_{r_{i-1}})\le  C \log i$ for any $i\ge 2$, hence $V(\xi_{r_i})\le \sum_{j=1}^i C \log j \le 2 C i \log i$ for all large $i$.

  Therefore  for any large $n$, there is a unique $j$ such that $r_j \le n < r_{j+1}$, and   $V(\xi_n)\le V(\xi_{r_j})+ C (n-r_j) \le 2 C j \log j  + C (r_{j+1}- r_j) \sim 4 C j \log j$. Since $\sqrt{n\log n} \ge \sqrt{r_j \log r_j} \sim   \sqrt{2} j \log j$, we see that for all large $n$, $V(\xi_n) \le 3 C \sqrt{n \log n}$. This yields $\Theta \le 3 C$ with positive probability. Since $\Theta$ is deterministic, we get   $\Theta \le 3C$, the desired upper bound.
  
It remains to show that  $({\mathcal Z}_i)$ can survive with positive probability. We will make a coupling between $(\#{\mathcal Z}_i)$ and a homogenous Galton-Watson process. Let $b_i:= r_i- r_{i-1}-1$ and $$\nu_i:= \sum_{|w|=b_i} 1_{\{  V(w_k)   \le C k , \forall \,  k \le b_i ,   V(w)   \ge 2 C \log i\}} 1_{\{\exists u\in \T_w: |u|= b_i+1,   V(u)   \le C \log i \}}.$$

Observe  that $(\#{\mathcal Z}_i)$ is stochastically larger than an inhomogeneous Galton-Watson process, say $(Z'_i)$, with $\nu_i$ as the reproduction law of the particles in the $(i-1)$th generation of $Z'$. It is enough to show that $(Z'_i)$ can survive with positive probability (if we choose $C$ sufficiently large). 

We claim that for any $K_0\ge 1$, \begin{equation} \label{K0C} \limsup_{C\to\infty} \limsup_{i\to\infty} \p(\nu_i\le K_0)  \le
q_\T ,\end{equation}

\noindent with $q_\T:=\p(\T \mbox{ is finite})\in [0, 1)$ is the extinction probability of the Galton--Watson tree $\T$. In fact, \eqref{K0C} is an equality, as $\nu_i=0$ for all large $i$ on the set of extinction of $\T$.

Assume for the moment \eqref{K0C}. Let $q'\in (q_\T, 1)$. Let $K_0$ be an integer such that $K_0 (1- q')>1$.    We may choose (and then fix) $C> 0$ large enough such that for all large $i$, $\p(\nu_i\le K_0)\le q'$.    Then for all large $i$, $\nu_i$ is stochastically larger that a Bernoulli variable which equals $K_0$ with probability $1-q'$ and $0$ with probability $q'$. Since the (homogenous) Galton-Watson process with such a Bernoulli variable as reproduction law is supercritical, we see that the inhomogeneous Galton-Watson process $(Z'_i)$ will survive with positive probability.

The rest of this subsection is devoted to the proof of \eqref{K0C}.

 Recall that $\e[S_1]={\tt m}>0$. Let $C> {\tt m}$ be large enough such that  \begin{equation} i^{C t^*} \, b_i^{-\alpha_2}\to \infty, \qquad \mbox{as } i\to\infty. \label{Ct*} \end{equation}

Let for $\ell\ge 0$, ${\mathcal F}_\ell:=\sigma\{u, V(u): |u|\le \ell\}$ be the $\sigma$-field generated by the BRW up to generation $\ell$. Conditionally on ${\mathcal F}_{b_i}$, $\nu_i$ is a sum of independent Bernoulli random variables: $$\e[\nu_i \, |\, {\mathcal F}_{b_i}] =\sum_{|w|=b_i} 1_{\{  V(w_k)   \le C k , \forall \,  k \le b_i ,   V(w)   \ge 2 C \log i\}}\,   \vartheta(V(w)-C \log i),$$ 

\noindent with $\vartheta(y):= \p(\exists   |u|=1 :  V(u)   \le - y)$. By assumption for all large $y$, $ \vartheta(y) \ge e^{-t^* y} y^{-\alpha_2}$. 
%
Therefore \begin{equation} \e[\nu_i \, |\, {\mathcal F}_{b_i}]\ge i^{C t^*} \, (C b_i)^{-\alpha_2} \, W'_i , \end{equation}

\noindent with $$W'_i:= \sum_{|w|=b_i} 1_{\{  V(w_k)   \le C k , \forall \,  k \le b_i ,   V(w)   \ge 2 C \log i\}} e^{-t^* V(w)}.$$

\noindent Observe that \begin{align*} \mbox{Var}(\nu_i \, |\, {\mathcal F}_{b_i}) 
&=  \sum_{|w|=b_i} 1_{\{  V(w_k)   \le C k , \forall \,  k \le b_i ,   V(w)   \ge 2 C \log i\}} \, \vartheta\left(V(w)-C \log i\right)\, \big(1-  \vartheta\left(V(w)-C \log i\right)\big)
\\
&\le 
  \e[\nu_i \, |\, {\mathcal F}_{b_i}].\end{align*}

\noindent This implies that for any $\varepsilon\in (0, 1)$, \begin{align*} 
\p\Big(\nu_i - \e[\nu_i \, |\, {\mathcal F}_{b_i}] \le - \varepsilon \e[\nu_i \, |\, {\mathcal F}_{b_i}], W'_i \ge \varepsilon\Big) 
&\le 
\e \Big[ \frac{\mbox{Var}(\nu_i \, |\, {\mathcal F}_{b_i})}{\varepsilon^2 (\e[\nu_i \, |\, {\mathcal F}_{b_i}])^2} 1_{\{W'_i \ge \varepsilon\}}\Big]
\\
&\le 
\e \Big[ \frac{1}{\varepsilon^2 \e[\nu_i \, |\, {\mathcal F}_{b_i}] } 1_{\{W'_i \ge \varepsilon\}}\Big]
\\
&\le  \frac{1}{\varepsilon^3 \,  i^{C t^*} (C b_i)^{-\alpha_2}}
\\
&\to    0  , 
\end{align*}

\noindent as $i\to \infty$, where the last convergence follows from \eqref{Ct*}.

   On the set $\{W'_i\ge \varepsilon\}$, $\e[\nu_i \, |\, {\mathcal F}_{b_i}] \ge \varepsilon i^{C t^*} (C b_i)^{-\alpha_2} \ge K_0/(1-\varepsilon)$ for all large $i$. Therefore,  $$\p(\nu_i\le K_0) 
\le
 \p\Big(\nu_i - \e[\nu_i \, |\, {\mathcal F}_{b_i}] \le - \varepsilon \e[\nu_i \, |\, {\mathcal F}_{b_i}], W'_i\ge \varepsilon\Big) + \p\Big( W'_i < \varepsilon\Big),$$
 
 \noindent yielding that \begin{equation}\limsup_{i\to\infty}\p(\nu_i\le K_0) 
 \le  \limsup_{i\to\infty} \p(W'_i<  \varepsilon). \end{equation}

To estimate $\p ( W'_i < \varepsilon )$, we compare $W'_i$ with $W_i$, where 
$$W_i:= \sum_{|w|=b_i}   e^{-t^* V(w)}.$$

\noindent Since by assumption  $ \e\left[ \sum_{|u|=1} e^{-t^* V(u)} \log_+ (\sum_{|u|=1} e^{-t^* V(u)})\right] < \infty  $, we may apply  Biggin's theorem (\cite{biggins-mart-cvg}) to see that  $W_i\to W_\infty$ which is positive $\p^*$-a.s.  By the many-to-one formula \eqref{many-to-one}, \begin{align*} \e[W_i-W'_i]&=  1- \p\Big( S_k \le C k, \forall\, k\le b_i, \, S_{b_i} \ge 2 C \log i\Big)
\\
&\le  \p\Big( \max_{k\ge 1} \frac{S_k}{k} > C\Big)+ \p\Big( S_{b_i} < 2 C \log i\Big) .\end{align*}

\noindent Recall that $\e[S_1]={\tt m} >0$ and $b_i \sim  2 i \log i$ as $i\to\infty$. The law of large numbers yields that $\p\big( S_{b_i} < 2 C \log i\big)\to0$. Using $\p(W'_i < \varepsilon) \le \p(W_i \le 2 \varepsilon) +  \p(W_i-W'_i>\varepsilon) \le
\p(W_i \le 2 \varepsilon) + \frac1\varepsilon  \e[W_i-W'_i]$, we get that   \begin{align*} \limsup_{i\to\infty} \p(\nu_i\le K_0) 
 &\le 
\limsup_{i\to\infty} \p(W'_i < \varepsilon)
\\
&\le 
\limsup_{i\to\infty} \p(W_i \le 2 \varepsilon) + \frac1\varepsilon \p\Big( \max_{k\ge 1} \frac{S_k}{k} > C\Big).
\end{align*}

\noindent Since $\p( \max_{k\ge 1} \frac{S_k}{k} > C) \to 0$ as $C\to\infty$ and $\p(W_\infty=0)=q_\T$, we complete the proof of \eqref{K0C} by letting $C\to\infty$ then $\varepsilon\to0$.  $\Box$

\appendix
\section{Estimates on random walks}\label{s:est}

\begin{lemma}\label{l:renewal1}
Let $(S_n,\, n\ge 0)$  be a real-valued random walk such that $\p(S_1\neq 0)>0$.   Then 
  $$
  \sup_{x\ge 0} \e\left[  \sum_{k=0}^\infty {\bf 1}_{\{\min_{0\le i\le k}S_i +x \ge 0\}} e^{- S_k-x}  \right]<\infty.
  $$
  \end{lemma}

\noindent {\it Proof}. 
Let $\upsilon_0:=0$ and $\upsilon_j:=\inf\{n> \upsilon_{j-1}: S_n  \le  S_{\upsilon_{j-1}}\}$ for $j\ge1$,  be the weak descending ladder times of $S$, with the convention that $\inf \emptyset=\infty$. We have 
$$
\e\left[  \sum_{k=0}^\infty {\bf 1}_{\{\min_{0\le i\le k} S_i +x \ge 0\}} e^{- S_k}  \right] = \sum_{j=0}^\infty \e\left[  {\bf 1}_{\{ S_{\upsilon_{j}}\ge -x,\, \upsilon_j<\infty \}} e^{-S_{\upsilon_j}}\right]\e\left[\sum_{k=0}^\infty e^{-S_k}{\bf 1}_{\{k<\upsilon_1\}}\right]. 
$$

\noindent The term $ \e\left[\sum_{k=0}^\infty e^{-S_k}{\bf 1}_{\{k<\upsilon_1\}}\right] $ can be evaluated through the ladder heights. Let $\upsilon^+_0:=0, \upsilon^+_j:= \inf\{i> \upsilon^+_{j-1}: S_i > S_{\upsilon^+_{j-1}}\}$ for $j\ge 1$. By the time-reversal for the random walk, $$ \e\left[\sum_{k=0}^\infty e^{-S_k}{\bf 1}_{\{k<\upsilon_1\}}\right]= \e \left[ \sum_{j=0}^\infty e^{-   S_{\upsilon^+_j}} {\bf 1}_{\{\upsilon^+_j<\infty\}}\right]= \sum_{j=0}^\infty  \left(\e \left[  e^{-   S_{\upsilon^+_1}} {\bf 1}_{\{\upsilon^+_1<\infty\}} \right]\right)^j <\infty$$ as $\e \left[  e^{- S_{\upsilon^+_1}} {\bf 1}_{\{\upsilon^+_1<\infty\}}\right]< 1$.  Finally,  $\sum_{j=0}^\infty \e\left[  {\bf 1}_{\{ S_{\upsilon_j}\ge -x,\,\upsilon_j<\infty \}} e^{-S_{\upsilon_j}}\right] \le \sum_{z=0}^{\lceil x \rceil } e^{z}U^-(z)$ where
$$
U^-(z):= \sum_{j=0}^\infty \p\left(S_{\upsilon_j} \in [-z,-z+1), \, \upsilon_j<\infty\right).
$$

\noindent  By the classical renewal theorem (see for example \cite[Theorem II.4.2]{Gut}), $\sup_{z\ge 0} U^-(z)=c<\infty$. Therefore
$$
\sum_{j=0}^\infty \e\left[  {\bf 1}_{\{ S_{\upsilon_j}\ge -x , \,\upsilon_j<\infty\}} e^{-S_{\upsilon_j}}\right] \le \frac{c \, e}{e-1} e^{\lceil x \rceil}.
$$

\noindent It completes the proof. $\Box$

\bigskip

  \begin{lemma}\label{l:hitting}
    Let $(S_n,\, n\ge 0)$  be a random walk with values in $\mathbb Z$ such that $\p(S_1\ge 2)=0$.  In agreement with \eqref{def:r}, let $\sigma_1:=\inf\{k\ge 1: S_k=0\}$. Suppose that $\e[S_1]=0$ and $\e[S_1^2]<\infty$. Then there exists $\theta>0$ such that   \begin{equation}    \p(T'_{-n}< \sigma_1)
  \sim \frac\theta{n}, \qquad n\to\infty, \label{eq:theta}\end{equation}
\noindent where $T'_{-n}:=\inf\{i\ge 0: S_i \le -n\}$.  Furthermore, there exists some sequence $ (K_n)_{n\ge 0}$ such that    $K_n\ge 1$, $K_n=o(n)$ as $n\to\infty$  and \begin{equation} \p(T'_{-n}< \sigma_1, \, |S_{T'_{-n}}+n| \le K_n )\sim \frac\theta{n}, \qquad n\to\infty . \label{eq:thetaKn}  \end{equation}
  \end{lemma}
  
  \noindent {\it Proof}. We will prove the lemma with $\theta:= \e \big[ |S_{\upsilon_1}|  1_{\{\upsilon_1 < \sigma_1\}}\big]$ where $\upsilon_1:=\inf\{i\ge 1: S_i < 0\}$.  The optional stopping theorem applied to  $T'_{-n}\land \sigma_1$ implies that for $-n< x <0$, 
  $$
 x=\e_x [S_{\sigma_1 \wedge T'_{-n}}]=\e_x[S_{T'_{-n}}{\bf 1}_{\{ T'_{-n} < \sigma_1\}} ] \le -n \p_x(T'_n <\sigma_1).
  $$
  
  \noindent  Hence $\p_x(T'_{-n} < \sigma_1)\le \frac{|x|}{n}$ for any $-n< x <0$, and also for any $x\le -n$. By the strong Markov property applied to $\upsilon_1$,
  \begin{equation}\label{eq: upper bound qn}
  \p\left(T'_{-n} < \sigma_1\right)= \e\left[{\bf 1}_{\{\upsilon_1 < \sigma_1\}} \p_{S_{\upsilon_1}}(T'_{-n} < \sigma_1)\right]  \le \frac{\theta}{n}.
  \end{equation}
  
\noindent     Discussing the last value taken by $S$ before leaving $[-(n-1),-1]$, we have for any $-n<x<0$, 
$$
\e_x[|S_{T'_{-n}}+n| {\bf 1}_{\{ T'_{-n}< \sigma_1\}}] = \sum_{k=1}^{n-1} G_n(x,-k) \e[|S_1+n-k|{\bf 1}_{\{S_1\le k-n\}}]
$$
  
  \noindent where $G_n(x,-k):=\sum_{\ell=0}^\infty \p_x(S_\ell=-k,\, \ell< T'_{-n}\land \sigma_1)$. Hence the inequality 
  \begin{equation}\label{eq:overshoot}
  \e_x[|S_{T'_{-n}}+n| {\bf 1}_{\{ T'_{-n}< \sigma_1\}}] \le \sum_{k=1}^{n-1} G_n(x,-k) \e[|S_1|{\bf 1}_{\{S_1\le k-n\}}].
  \end{equation}
  
\noindent  As in \eqref{TA}, let $T_{\{-k\}}=\inf\{i\ge 1\,:\, S_i=-k\}$. In the sequel, the constant $c$ does not depend on $n,k,x$ and can change its value from line to line. By \eqref{eq:ruin}, for all  $-k<x<0$
  $$
  \p_x(T_{\{-k\}}<   \sigma_1) \le c \frac{|x|}{k}.
  $$
  
  \noindent We can suppose that the inequality stays true  when $x\le -k$ by taking a bigger $c$ if necessary. Moreover,  by \eqref{eq:ruin},  for all $0<k<n$,
  $$
  \p_{-k}( T_{\{-k\}} > T'_{-n} \land \sigma_1) \ge c \max\big(\frac1k,\frac1{n-k}\big).
  $$

\noindent   The last two inequalities imply that $G_n(x,-k)\le c \min(k,n-k) \frac{|x|}{k}$ for all $-n<x<0$ and $0<k<n$. Splitting the sum on the right-hand side of \eqref{eq:overshoot} according to whether $k\le \lfloor \frac{n}{2} \rfloor$ or $k> \lfloor \frac{n}{2} \rfloor$ gives
\begin{equation}   
\e_x[|S_{T'_{-n}}+n| {\bf 1}_{\{ T'_{-n}< \sigma_1\}}] \le c |x|  (I_1+I_2) \label{I1I2}
 \end{equation}
  
  \noindent where
  $$
I_1:= \sum_{k=1}^{\lfloor \frac{n}{2} \rfloor} \e[|S_1|{\bf 1}_{\{S_1\le k-n\}}],  \qquad I_2:= \frac1n \sum_{k=1}^{\lfloor \frac{n}{2} \rfloor} k\, \e[|S_1|{\bf 1}_{\{S_1\le -k\}}].
$$

\noindent We observe that 
$$
I_1 \le \lfloor \frac{n}{2} \rfloor  \e[|S_1|{\bf 1}_{\{2S_1\le -n\}}]  \le \e[S_1^2 {\bf 1}_{\{2S_1\le -n\}}] \to 0
$$

\noindent as $n\to\infty$ by dominated convergence. Similarly, since $\sum_{k=1}^{\min(p,q)} k \le c p\min(p,q)$ for $p=|S_1|$ and $q=\lfloor \frac{n}{2}\rfloor $, 
$$
I_2 \le c\,  \e\bigg[S_1^2 \frac{\min(|S_1|,n)}{n}\bigg] \to 0.
$$

\noindent   We deduce that $\lim_{n\to\infty} \e_x[|S_{T'_{-n}}+n|\, {\bf 1}_{\{ T'_{-n}< \sigma_1\}}] =0$. Writing
 $$
 x=\e_x [S_{\sigma_1 \wedge T'_{-n}}]=-n\p_x(T'_{-n}<\sigma_1) + \e_x[(S_{T'_{-n}}+n)\, {\bf 1}_{\{ T'_{-n}< \sigma_1\}}],
 $$

\noindent  it implies that for all fixed $x<0$, $\p_x(T'_{-n} < \sigma_1)\sim \frac{|x|}{n}$.  By Fatou's lemma,
$$
\liminf_{n\to\infty}  n\p\left(T'_{-n} < \sigma_1\right)= \liminf_{n\to\infty} \e\left[{\bf 1}_{\{\upsilon_1 < \sigma_1\}} n \p_{S_{\upsilon_1}}(T'_{-n} < \sigma_1)\right]  \ge  \theta.
$$

\noindent Together with \eqref{eq: upper bound qn}, it completes the proof of \eqref{eq:theta}.

To prove \eqref{eq:thetaKn}, it is sufficient to show that we can find a sequence $K_n=o(n)$ with $K_n\ge 1$, such that 
\begin{equation}  \p(T'_{-n}< \sigma_1, \, |S_{T'_{-n}}+n| > K_n )  =o\left(\frac1{n}\right), \qquad n\to\infty.  \label{proba>Kn} \end{equation}

\noindent In fact, the probability term in \eqref{proba>Kn} is less than \begin{align*}     \p\left(S_{\upsilon_1} \le -n\right)+  \e\left[{\bf 1}_{\{\upsilon_1 < \sigma_1, S_{\upsilon_1} > -n \}} \p_{S_{\upsilon_1}}(T'_{-n} < \sigma_1, \, |S_{T'_{-n}}+n| > K_n)\right] 
=: I_3+ I_4.
\end{align*}

\noindent By the Markov inequality, $$I_3\le  \frac1{n}\, \e \left[ |S_{\upsilon_1}| {\bf 1}_{\{S_{\upsilon_1} \le -n\}}\right].$$

\noindent According to Doney  \cite[Corollary 1]{Doney}, $\e \left[ |S_{\upsilon_1}|  \right]< \infty$, hence $$I_3= o\left(\frac1{n}\right), \qquad n\to \infty.$$

For $I_4$, using \eqref{I1I2}, we have \begin{align*}   I_4 &\le
 \frac1{K_n} \e\left[{\bf 1}_{\{\upsilon_1 < \sigma_1, S_{\upsilon_1} > -n \}}  \e_{S_{\upsilon_1}}\big(|S_{T'_{-n}}+n| {\bf 1}_{\{ T'_{-n}< \sigma_1\}} \big)\right]
 \\
 &\le \frac{c(I_1+I_2)}{K_n} \e\left[ |S_{\upsilon_1}|\right]
  \le \frac{c' \varepsilon_n}{K_n},
 \end{align*} 
 
 \noindent with some positive constant $c'$ and $\varepsilon_n:= \e[S_1^2 {\bf 1}_{\{2S_1\le -n\}}]+\e\left[S_1^2 \frac{\min(|S_1|,n)}{n}\right]  \to0.$
 
 If we choose $K_n:= \max(1, n \sqrt{\varepsilon_n})$, then $K_n=o(n)$ and $I_4 \le \frac{c'}{n} \varepsilon_n^{1/2} =o\left(\frac1{n}\right)$.  Thus,  we obtain \eqref{proba>Kn} and complete the proof of Lemma \ref{l:hitting}.
$\Box$


\end{document}